\documentclass[onefignum,onetabnum]{siamart220329}
\usepackage{bm}
\usepackage{algorithmic}
\usepackage{algorithm}
\usepackage{mathrsfs}
\usepackage{epsfig}
\usepackage{enumitem}
\usepackage{diagbox}
\usepackage{cases}

\usepackage{booktabs}
\usepackage{subfigure}
\usepackage{amssymb}
\usepackage{dsfont}

\usepackage{stmaryrd}

\usepackage{xparse}

\NewDocumentCommand{\dgal}{sO{}m}{%
  \IfBooleanTF{#1}
    {\dgalext{#3}}
    {\dgalx[#2]{#3}}%
}

\NewDocumentCommand{\dgalext}{m}{%
  \sbox0{%
    \mathsurround=0pt 
    $\left\{\vphantom{#1}\right.\kern-\nulldelimiterspace$%
  }%
  \sbox2{\{}%
  \ifdim\ht0=\ht2
    \{\kern-.45\wd2 \{#1\}\kern-.45\wd2 \}%
  \else
  \fi
}

\NewDocumentCommand{\dgalx}{om}{%
  \sbox0{\mathsurround=0pt$#1\{$}%
  \sbox2{\{}%
  \ifdim\ht0=\ht2
    \{\kern-.45\wd2 \{#2\}\kern-.45\wd2 \}%
  \else
    \mathopen{#1\{\kern-.5\wd0 #1\{}
    #2
    \mathclose{#1\}\kern-.5\wd0 #1\}}
  \fi
}

\usepackage{lipsum}
\usepackage{amsfonts}
\usepackage{graphicx}
\usepackage{epstopdf}
\usepackage{algorithmic}

\usepackage{animate}
\usepackage{amsmath}
\usepackage{amssymb}
\usepackage{graphics}
\usepackage{graphicx}
\usepackage{footnote}
\usepackage{textcomp}
\usepackage{mathrsfs}
\usepackage{epstopdf}
\usepackage{array}
\usepackage[maxfloats=99]{morefloats}
\usepackage{url}
\usepackage{cases}
\usepackage{mathscinet}
\usepackage[normalem]{ulem}
\usepackage{algorithmic, algorithm}
\usepackage{color}
\usepackage{cite}

\ifpdf
  \DeclareGraphicsExtensions{.eps,.pdf,.png,.jpg}
\else
  \DeclareGraphicsExtensions{.eps}
\fi

\newsiamremark{remark}{Remark}
\newsiamremark{hypothesis}{Hypothesis}
\crefname{hypothesis}{Hypothesis}{Hypotheses}
\newsiamthm{claim}{Claim}

\headers{Temporal Difference Methods}{}

\title{Temporal Difference Learning for High-Dimensional PIDEs with Jumps}

\author{Liwei Lu\thanks{Yau Mathematical Sciences Center, Tsinghua University, Beijing, 100084, China (\email{llw20@mails.tsinghua.edu.cn}).}
\and 
Hailong Guo\thanks{ School of Mathematics and Statistics,  The University of Melbourne,  Parkville, VIC 3010, Australia (\email{hailong.guo@unimelb.edu.au}).}
\and 
Xu Yang\thanks{Department of Mathematics, University of California, Santa Barbara, CA 93106, USA (\email{xuyang@math.ucsb.edu}).}
\and 
Yi Zhu\thanks{Corresponding author, Yau Mathematical Sciences Center, Tsinghua University, Beijing, 100084, China and Beijing Institute of Mathematical Sciences and Applications, Beijing, 101408, China (\email{yizhu@tsinghua.edu.cn}).}
}

\usepackage{amsopn}

\ifpdf
\hypersetup{
  pdftitle={Temporal Difference Learning for PIDEs},
  pdfauthor={Liwei Lu, Hailong Guo, Xu Yang, Yi Zhu}
}
\fi

\usepackage{multirow}

\graphicspath{{fig/}{fig/mesh/}}

\usepackage{stmaryrd}

\newcommand{\R}{\mathbb{R}}

\newcommand{\E}{\mathbb{E}}

\newcommand{\dif}{\mathop{}\!\mathrm{d}}

\begin{document}

\maketitle

\begin{abstract}
  In this paper, we propose a deep learning framework for solving high-dimensional partial integro-differential equations (PIDEs) based on the temporal difference learning. We introduce a set of L\'{e}vy processes and construct a corresponding reinforcement learning model. To simulate the entire process, we use deep neural networks to represent the solutions and non-local terms of the equations. Subsequently, we train the networks using the temporal difference error, termination condition, and properties of the non-local terms as the loss function. The relative error of the method reaches $\mathcal{O}(10^{-3})$ in 100-dimensional experiments and $\mathcal{O}(10^{-4})$ in one-dimensional pure jump problems. Additionally, our method demonstrates the advantages of low computational cost and robustness, making it well-suited for addressing problems with different forms and intensities of jumps.
\end{abstract}




\section{Introduction}

Partial integro-differential equation (PIDE) has found widespread applications in various real-world scenarios, including engineering, biology, and particularly finance\cite{abergel2010nonlinear,kavallaris2018non,goswami2016system}. Examples of such phenomena can be observed in Ohmic Heating, Resistance Spot Welding, Gierer-Meinhardt System, and option pricing\cite{kavallaris2018non,goswami2016system}. Classical methods such as the finite difference method\cite{cont2005finite,kwon2011second} and finite elements method\cite{pablo2023finite,goswami2013optimal} have been traditionally used to solve PIDEs, and they exhibit good performance in low-dimensional problems. However, these mesh-based approaches face significant challenges in high-dimensional cases due to the curse of dimensionality\cite{koppen2000curse}. The computational complexity grows exponentially as the number of dimensions increases, rendering these methods impractical. Nevertheless, many real-world models formulated as PIDEs inherently possess high-dimensional characteristics. For instance, in option pricing, the dimensions correspond to the number of underlying assets, while in reaction-diffusion systems of cell dynamics, they correspond to the number of reacting substances. 

Neural networks have emerged as a promising tool for alleviating the curse of dimensionality in solving high-dimensional partial differential equations (PDEs) without non-local terms. In recent years, several neural network-based methods have been widely adopted for this purpose\cite{cai2021physics,cuomo2022scientific}. The Physics-Informed Neural Network (PINN) \cite{raissi2019physics}\cite{chiu2022can,kharazmi2021hp,li2022revisiting,mcclenny2023self,moseley2021finite,wight2020solving,yu2022gradient,karniadakis2021physics,lu2021deepxde} incorporates the residual term of the PDE into the training process and utilize automatic differentiation to handle the differential operations, and it has sparked a wave of using deep learning to solve PDEs. The idea of Deep Galerkin Method (DGM) \cite{sirignano2018dgm} is similar to PINN, but it incorporates the norm of the residual term in its loss function and adopts $L_2$ norm. The Deep Mixed Residual Method (MIM) \cite{lyu2022mim} rewrites high-order PDEs into a system of first-order PDEs for solving purposes. The Weak Adversarial Network (WAN) \cite{zang2020weak} leverages the analogy between the weak form of the PDE and generative adversarial network, integration by parts to transfer derivatives onto the test functions. The Deep Ritz Method \cite{yu2018deep} utilizes the Ritz variational form to solve high-dimensional PDEs and eigenvalue problems. The Deep Nitsche method \cite{liao2019deep} represents the trial functions by neural network and employs the Nitsche variational form to handle mixed boundary value problems. In addition, introducing stochastic differential equations (SDEs) related to PDE provides another potential avenue for mitigating the curse of dimensionality. Deep Backward Stochastic Differential Equation (Deep BSDE) \cite{han2017deep,han2018solving} approximates the gradients of the solution at different time points using a series of neural networks and trains them with the terminal condition. The Forward-Backward Stochastic Neural Network (FBSNN) \cite{raissi2018forward} incorporates trajectory information into the training process to directly obtain the solution of the equation. \cite{zeng2022deep} involves using two neural networks to approximate the solution and its gradient, or by auto-differentiation, and performing training at each time point based on temporal difference method.

When it comes to partial integro-differential equations, the literature on deep neural network-based approaches is not as extensive. \cite{gnoatto2022deep} extends the work of \cite{han2017deep,han2018solving} to a more general case involving jump processes and uses additional neural networks to approximate non-local terms. For problems with a finite number of jumps, \cite{frey2022deep} circumvents the computation of non-local terms by introducing a stochastic process with integration on a Poisson random measure. In \cite{castro2022deep}, an error estimation technique is proposed for approximating forward-backward stochastic systems using neural networks. Furthermore, \cite{wang2023deep} approximates the initial value function, gradients of the solution, and integral kernel using neural networks, with the termination condition serving as a constraint.

This paper aims to numerically compute solutions to partial integro-differential equations. A group of L\'{e}vy-type forward-backward stochastic processes is introduced based on the target PIDE. Inspired by \cite{zeng2022deep}, a reinforcement learning framework is established. The equation's solution and non-local terms are represented using neural networks. Subsequently, a loss function is constructed to update the network parameters, taking into account the errors from temporal difference methods, termination conditions, and the properties of non-local terms. This method exhibits two primary advantages: 1) It exhibits a reduced computational cost. Capitalizing on the benefits of temporal difference learning, the proposed approach eliminates the need to wait for the completion of an entire trajectory simulation before updating parameters. Moreover, as the dimensionality increases, the required number of trajectories does not experience significant growth. 2) The method has fast convergence and high precision. With rapid convergence, the error remains within the range of $\mathcal{O}(10^{-4})$ for one-dimensional pure jump problems and $\mathcal{O}(10^{-3})$ for high-dimensional problems.

The rest of this paper is organized as follows. Section \ref{sec.method} presents the overall methodology, which involves the incorporation of L\'{e}vy processes and reinforcement learning models to establish the training objectives. Section \ref{sec.numerical} demonstrates numerical examples for solving one-dimensional and high-dimensional PIDEs. Finally, in Section \ref{sec.conclusion}, we provide a summary of the entire approach.


\section{Methodology}  \label{sec.method}

The Brownian motion, defined based on the Gaussian distribution, has been widely adopted for modeling noise. However, in the scenarios of real applications, non-Gaussian noises are better characterized by L\'{e}vy motion, which is particularly relevant in fields like finance, chemistry, engineering, and geophysics\cite{delong2013backward,duan2015introduction}. L\'{e}vy processes are stochastic processes that satisfy zero initial value, stationary independent increments, and stochastically continuous sample paths\cite{applebaum2009levy}. For any L\'{e}vy process $L_t$, the Poisson random measure $N$, the compensated Poisson random measure $\widetilde{N}$ and the jump measure $\nu$ can be defined by
\begin{gather}
  N(t,S)(\omega) := \mbox{card}\{s\in[0,t): L_s(\omega)-L_{s-}(\omega)\in S]\}, \\ 
  \nu(S) := \E N(1, S)(\omega), \\ 
  \widetilde{N}(t,S) = N(t,S) - t\nu(S),
\end{gather}
where $S$ is a Borel set on $\R^d\backslash\{0\}$.

The L\'{e}vy form of the It\^{o}'s formula bridges the gap between integro-differential operators and L\'{e}vy processes. Nevertheless, there remains a considerable disparity in simulating the L\'{e}vy process and accurately computing the non-local terms in the operator. Furthermore, efficiently solving high-dimensional problems at a lower cost presents a formidable challenge. To tackle these issues, we propose a framework that employs reinforcement learning and utilizes deep neural networks.

\subsection{PIDE and L\'{e}vy process}
In the present investigation, our attention is directed towards the resolution of the following partial integro-differential equation (PIDE)
\begin{equation} \label{eq.PIDE}
  \left\{
  \begin{aligned}
      & \frac{\partial u}{\partial t} + b\cdot\nabla u + \frac{1}{2}\mbox{Tr}(\sigma\sigma^T\mbox{H}(u)) + \mathcal{A}u + f = 0, \\ 
      & u(T, \cdot) = g(\cdot),
  \end{aligned}
  \right.
\end{equation}
where $ \nabla u $ and $ \mbox{H}(u) $ correspond to the gradient and Hessian matrix of the function $u(t,x)$ with respect to the spatial variable $x\in\R^d$, and $0<t<T$. The vectors $ b=b(x)\in\R^d $, matrix $ \sigma=\sigma(x)\in\R^{d\times d} $ as well as the function $ f=f(t,x,u,\sigma^T\nabla u)\in\R $, are all pre-defined. Here
\begin{equation}
  \mathcal{A}u(t,x) = \int_{\R^d}(u(t,x+G(x,z)) - u(t,x) - G(x,z)\cdot\nabla u(t,x))\nu(\dif z),
\end{equation}
where $ G=G(x,z)\in\R^d\times\R^d\to\R^d $, and $ \nu $ is a L\'{e}vy measure associated with a Poisson random measure $ N $.

Traditional approaches, such as finite element and finite difference methods, tend to be insufficient in dealing with non-local terms and high-dimensional scenarios. To overcome these limitations, we propose the following L\'{e}vy-type process on the probability space $(\Omega, \mathcal{F}, \mathbb{P})$ in the It\^{o} sense
\begin{equation} \label{eq.Xt}
  \dif X_t = b(X_{t})\dif t + \sigma(X_{t})\dif W_t + \int_{\R^d} G(X_{t},z) \widetilde{N}(\dif t,\dif z).
\end{equation}
In this context, $\{W_t\}_{t=0}^T$ denotes a $d$-dimensional Brownian motion, and $ \widetilde{N}(\dif t,\dif z)=N(\dif t,\dif z) - \nu(\dif z)\dif t $ represents the compensated Poisson measure. $\{\mathcal{F}_t\}_{t=0}^T$ is the filtration generated by Brownian motion $W_t$ and Poisson random measure $ N $. If $X_0$ is assumed to follow a known single-point distribution $ X_0=\xi\in\R^d $, the entire process $X_t$ can be completely known and simulated. Subsequently, we introduce the following three stochastic processes:
\begin{equation} \label{eq.YZU}
  Y_t = u(t,X_t), \ Z_t = \nabla u(t,X_t), \ U_t = \int_{\R^d}(u(t,X_t+G(X_t,z)) - u(t,X_t))\nu(\dif z).
\end{equation}
According to It\^{o} formula for L\'{e}vy-type stochastic integrals\cite{applebaum2009levy},
\begin{equation} \label{eq.Yt}
    \begin{aligned}
        \dif Y_t = & \left[\frac{\partial u}{\partial t}(t,X_t)+b(X_t)\cdot\nabla u(t,X_t) \right. \\
        & +\left.\frac{1}{2}\mbox{Tr}\left(\sigma(X_t)\sigma(X_t)^T\mbox{H}(u)(t,X_t)\right) + \mathcal{A}u(t,X_t) \right]\dif t \\ 
        & + \nabla u(t,X_t)\cdot\sigma(X_t) \dif W_t + \int_{R^d} \left[u(t, X_t+G(X_t,z)) - u(t,X_t)\right] \widetilde{N}(\dif t,\dif z) \\ 
        = & -f(t,X_t,Y_t,\sigma(X_t)^T Z_t)\dif t + (\sigma(X_t)^T Z_t)^T\dif W_t \\ 
        & + \int_{R^d} \left[u(t, X_t+G(X_t,z)) - u(t,X_t)\right] \widetilde{N}(\dif t,\dif z), \\ 
    \end{aligned}
\end{equation}
with $Y_T=g(X_T)$.

Equations \eqref{eq.Xt} and \eqref{eq.Yt} are collectively referred to as the forward-backward stochastic differential equation (FBSDE) system\cite{raissi2018forward}, which establishes a fundamental link between solving the equation and simulating the stochastic process. In other words, if $u(t,x)$ represents a solution to equation \eqref{eq.PIDE}, the processes $(X_t, Y_t, Z_t, U_t)$ defined through \eqref{eq.Xt} and \eqref{eq.YZU} are required to satisfy equation \eqref{eq.Yt}. This provides a foundation for assessing the accuracy of the approximate solution $u(t,x)$.

In order to simulate these stochastic processes, the time $[0,T]$ is divided into $N$ equal subintervals: $ 0=t_0<t_1<t_2<\cdots<t_N=T $, resulting in the following numerical discretization of \eqref{eq.Xt} 
\begin{equation} \label{eq.dXt}
  \begin{gathered}
    X_{t_n+1} = X_{t_n} + b(X_{t_n})\Delta t + \sigma(X_{t_n})\Delta W_n \\
    + \sum_{i=N_n+1}^{N_{n+1}} G(X_{t_n},z_i) - \Delta t\int_{\R^d} G(X_{t_n},z) \nu(\dif z),
  \end{gathered}
\end{equation}
and discretization of \eqref{eq.Yt}
\begin{equation} \label{eq.dYt}
  \begin{gathered}
    Y_{t_{n+1}} = Y_{t_n}-f(t_n,X_{t_n},Y_{t_n},\sigma^T Z_{t_n})\Delta t + (\sigma(t_n,X_{t_n})^T Z_{t_n})^T\Delta W_n \\ 
    + \sum_{i=N_n+1}^{N_{n+1}} [u(t_n, X_{t_n}+G(X_{t_n},z_i))-u(t_n,X_{t_n})] - \Delta t U_{t_n}.
  \end{gathered}
\end{equation}
The terms $\Delta W_n$ in both \eqref{eq.dXt} and \eqref{eq.dYt} correspond to a same sample drawn from normal distribution $\mathcal{N}(0, \Delta t)$. $N_n$ represents the number of jumps that occur in the interval $[0, t_n]$ for the Poisson random measure $N$, while $z_i\in\R^d$ denotes the size of the $i$-th jump if arranging all jumps in the interval $[0, T]$ in chronological order. This discretization provides a practical framework for numerical approximation.

\subsection{Neural networks approximations} \label{sec.neural}

The discretized versions of equations \eqref{eq.dXt} and \eqref{eq.dYt} offer a viable approach for simulating the stochastic processes. However, several aspects remain unclear and require further investigation before achieving the final solution. While the simulation of the stochastic process $X_t$ can proceed smoothly once the initial value $X_0=\xi$ is determined, the simulation of $Y_t$ presents greater challenges due to its dependence on $Y_t$ itself, as well as $Z_t$ and $U_t$. The computation of the integral form in $U_t$ is a matter that needs careful consideration, and accurately calculating the difference in $Z_t$ is also not a trivial task.

In this work, neural networks have been adopted as an approximation for the solution $u(t,x)$ of the PIDE \eqref{eq.PIDE}. Specifically, when the values of $X_t$ are known, they can be readily fed into the neural network to obtain an approximation of $Y_t$. Moreover, the automatic differentiation technique enables us to conveniently approximate $Z_t$. However, the computation of $U_t$ remains a challenging task. Employing Monte Carlo methods to evaluate the integral in $U_t$ might not be an optimal choice, as it necessitates repetitive calculations for every possible value of $X_t$ and would require re-computation whenever $u(t,x)$ undergoes changes.

\begin{figure}[htbp]
  \centering
  \includegraphics[width=37em]{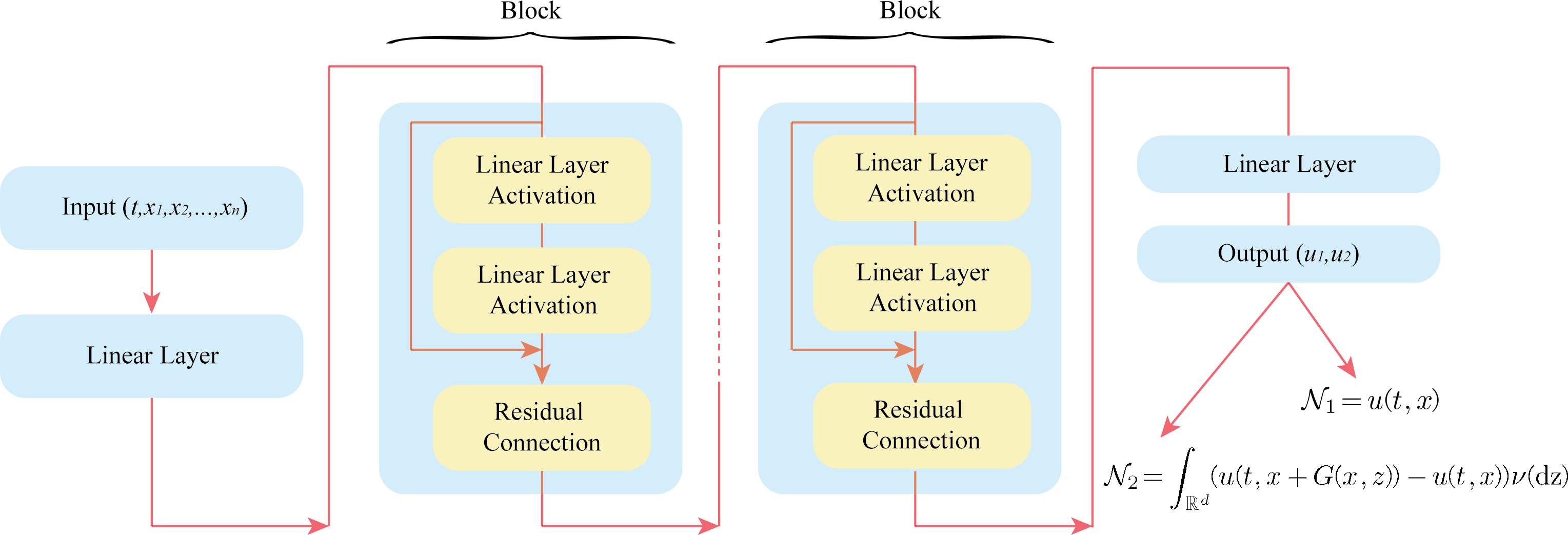}
  \caption{Architecture of the residual network utilized in this work. The neural network takes a $d+1$-dimensional input consisting of time $t\in\R$ and spatial variable $x\in\R^d$. The two outputs represent $u(t,x)$ and $\int_{\R^d}(u(t,x+G(x,z)) - u(t,x))\nu(\dif z)$ respectively.} \label{fig.net}
\end{figure}

To address the challenges associated with the computation of $U_t$, a modification is made to the output of the neural network. Specifically, two components are introduced: one output $\mathcal{N}_1$ to represent $u(t,x)$ and the other output $\mathcal{N}_2$ to represent $\int_{\R^d}(u(t,x+G(x,z)) - u(t,x))\nu(\dif z)$. It is important to note that the value of $\mathcal{N}_2$ depends not only on the input $(t, x)$ but also on the current $u(t,x)$. Consequently, employing a neural network with two outputs may be more suitable than using two separate networks, as it effectively captures the inherent relationship between these outputs. The structure of the neural network is depicted in Figure \ref{fig.net}, and it can be mathematically described as  
\begin{equation}
  \left\{
  \begin{aligned}
    & \bm{a}^1 = \bm{W}^1\bm{a}^0 +\bm{b}^1, \quad \bm{a}^0=(t,x)\in\R^{d+1}, \\ 
    & \bm{a}^i = h(\bm{W}^{2i-1}h(\bm{W}^{2i-2}\bm{a}^{i-1} + \bm{b}^{2i-2}) +\bm{b}^{2i-1}) + \bm{a}^{i-1}, \quad i=2,3,\cdots,l \\ 
    & \bm{a}^{l+1} = \bm{W}^{2l}\bm{a}^l +\bm{b}^{2l} \in \R^2.
  \end{aligned}
  \right. .
\end{equation}
Here, the $d+1$-dimensional input $(t, x)$ undergoes a linear layer, several residual blocks, and another linear layer, resulting in a 2-dimensional output. Each residual block consists of two linear layers combined with an activation function $h$, along with a residual connection.

Upon utilizing a neural network as an approximation, the following equation 
\begin{equation}
  \begin{gathered}
      Y_{t_{n+1}} = Y_{t_n}-f(t_n,X_{t_n},\mathcal{N}_1(t_n,X_n),\sigma^T \nabla\mathcal{N}_1(t_n,X_n))\Delta t \\
      + (\sigma(t_n,X_{t_n})^T \nabla\mathcal{N}_1(t_n,X_n))^T\Delta W_n \\ 
      + \sum_{i=N_n+1}^{N_{n+1}} [\mathcal{N}_1(t_n, X_{t_n}+G(X_{t_n},z_i))-\mathcal{N}_1(t_n,X_{t_n})] - \Delta t \mathcal{N}_2(t_n, X_{t_n})
  \end{gathered}
\end{equation}
allows us to simulate $Y_t$, where the definition of $z_i$ and $N_n$ remains consistent with equation \eqref{eq.dYt}. Consequently, it becomes necessary to impose a constraint on $\mathcal{N}_2$. There are various methods to enforce this constraint, and in this study, it is observed that 
\begin{equation}
  \left\{ \int_{0}^t \int_{\R^d}(u(t,X_t+G(X_t,z)) - u(t,X_t))\widetilde{N}(\dif t,\dif z) \right\}_{t=0}^T
\end{equation}
behaves as a martingale, implying
\begin{equation}
  \E \left[ \int_{t_n}^{t_{n+1}} \int_{\R^d} u(t,X_t + G(X_t,z)) - u(t,X_t) \widetilde{N}(\dif t,\dif z) \bigg| \mathcal{F}_{t_n} \right] = 0,
\end{equation}
within any time $[t_n, t_{n+1}]\subset[0,T]$. Inspired from this observation, we incorporate
\begin{equation} \label{eq.N2loss}
  \left|
  \frac{1}{M} \sum_{j=1}^M\left(\sum_{i=N_n+1}^{N_{n+1}} [\mathcal{N}_1(t_n, x^j_{t_n}+G(x^j_{t_n},z^j_i))-\mathcal{N}_1(t_n,x^j_{t_n})] - \Delta t \mathcal{N}_2(t_n, x^j_{t_n})\right)
  \right|
\end{equation}
into the final loss function, where $M$ denotes the number of training samples, and $x^j_{t_n}$ refers to the $j$-th sample of $X_{t_n}$. The constraint expressed in \eqref{eq.N2loss} establishes a link between $\mathcal{N}_1$ and $\mathcal{N}_2$. Consequently, our subsequent objective is to identify a criterion that enables $\mathcal{N}_1$ to closely approximate the solution $u(t,x)$ of the PIDE \eqref{eq.PIDE}.

\subsection{Temporal difference learning}

As an approach for learning through interaction, reinforcement learning (RL) has gained increasing attention in various fields, including solving high-dimensional PDEs\cite{liang2022finite}. Temporal difference (TD) is a reinforcement learning method with form of bootstrapping. It offers the advantage of not requiring prior knowledge about the environment and allows for learning from each state transition without waiting for the completion of a full trajectory. Remarkably, equation \eqref{eq.dYt} delineates the process $Y_t$ in an incremental manner, aligning harmoniously with the fundamental principles of temporal difference learning.

This observation leads to the establishment of a reinforcement learning model for our problem by defining a triplet $(S, R, P)$. The state $S=[0,T]\times\R^d$ encompasses the set of all possible spatiotemporal pairs $s=(t,x)$ under the initial condition $X_0=\xi$. The transition probability 
\begin{equation}
  P_{s_1s_2} = \mathbb{P}(S_{2}=(t_2,x_2)|S_{1}=(t_1,x_1)) = \mathbb{P}(X_{t_2}=x_2|X_{t_1}=x_1)
\end{equation}
represents the probability of transitioning from $x_1$ at time $t_1$ to $x_2$ at time $t_2$. Since the state process $X_t$ fully determines the transition probability $P$, the specific values of $P$ can be calculated using equation \eqref{eq.dXt}, making it a model-based problem. However, we will assume that this remains a model-free problem as our approach does not require knowledge of the transition probability $P$. Thus, our method can be generalized to more scenarios.

Furthermore, we aim to define a reward $R$ such that the corresponding state value function of the model aligns precisely with the solution $u(t,x)$ of PIDE \eqref{eq.PIDE}. This allows us to transform the solution of \eqref{eq.PIDE} into the computation of the value function within the framework of reinforcement learning. If $u(t,x)$ is employed as the value function in RL model, the process $Y_t$ defined in \eqref{eq.dYt} corresponds to the cumulative reward associated with the state process $X_t$. Hence, the reward can be defined by 
\begin{equation}
  \begin{gathered}
    R_{s_{n}s_{n+1}}= -f(t_n,X_{t_n},Y_{t_n},\sigma^T Z_{t_n})\Delta t + (\sigma(X_{t_n})^T Z_{t_n})^T\Delta W_n \\ 
    + \sum_{i=N_n+1}^{N_{n+1}} [u(t_n, X_{t_n}+G(X_{t_n},z_i))-u(t_n,X_{t_n})] - \Delta t U_{t_n},
  \end{gathered}
\end{equation}
using the incremental formulation in \eqref{eq.dYt} for $Y_t$. And in the case of neural network approximation, it can be computed by
\begin{equation}
  \begin{gathered}
    R_{s_ns_{n+1}} = -f(t_n,X_{t_n},\mathcal{N}_1(t_n,X_n),\sigma^T \nabla\mathcal{N}_1(t_n,X_n))\Delta t \\
    + (\sigma(X_{t_n})^T \nabla\mathcal{N}_1(t_n,X_n))^T\Delta W_n \\ 
    + \sum_{i=N_n+1}^{N_{n+1}} [\mathcal{N}_1(t_n, X_{t_n}+G(X_{t_n},z_i))-\mathcal{N}_1(t_n,X_{t_n})] - \Delta t \mathcal{N}_2(t_n, X_{t_n}).
  \end{gathered}
\end{equation}

Temporal difference learning updates the current value function at each time step without requiring the completion of the entire trajectory over the time interval $[0, T]$. Therefore, in this problem, the loss function is computed and optimization is performed at each time step. Since the current model does not involve actions and policies, equation 
\begin{equation} \label{eq.TD}
  u(s_{n}) \leftarrow u(s_{n}) + \alpha(R_{s_ns_{n+1}} + u(s_{n}) - u(s_{n+1}))
\end{equation}
is utilized to update the value function $u(t,x)$ in classical temporal difference method. The TD error $R_{s_ns_{n+1}} + u(s_{n}) - u(s_{n+1})$ serves as a measure of the accuracy of the current value function, which can be effectively used to impose constraints on the neural network $\mathcal{N}_1$. Consequently, the TD error is incorporated as the first component of the loss function
\begin{gather}
  \mbox{TD\_Error}^j_{t_n} = -f(t_n,x^j_{t_n},\mathcal{N}_1(t_n,x^j_n),\sigma^T \nabla\mathcal{N}_1(t_n,x^j_n))\Delta t \label{eq.TDerror} \\ 
  + (\sigma(x^j_{t_n})^T \nabla\mathcal{N}_1(t_n,x^j_n))^T\Delta W_n
  + \sum_{i=N_n+1}^{N_{n+1}} [\mathcal{N}_1(t_n, x^j_{t_n}+G(x^j_{t_n},z^j_i))-\mathcal{N}_1(t_n,x^j_{t_n})] \nonumber \\ 
  - \Delta t \mathcal{N}_2(t_n, x^j_{t_n})
  + \mathcal{N}_1(t_{n},x^j_{t_{n}}) - \mathcal{N}_1(t_{n+1},x^j_{t_{n+1}}), \nonumber \\ 
  \mbox{Loss}^1_{t_n} = \frac{1}{M} \sum_{j=1}^M |\mbox{TD\_Error}^j_{t_n}|^2. \label{eq.loss1}
\end{gather}
Here, $M$ represents the number of samples, and \eqref{eq.TDerror} represents the TD error of the $j$-th sample at time $t_n$.

Next, since the cumulative reward process $Y_t$ needs to satisfy not only equation \eqref{eq.dYt} but also $Y_T=g(X_T)$, the termination conditions and their gradients are incorporated as the second and third components of the loss function respectively, and evenly distributed across each time step $ 0=t_0<t_1<t_2<\cdots<t_N=T$,
\begin{gather}
  \mbox{Loss}^2_{t_n} = \frac{1}{N}\frac{1}{M} \sum_{j=1}^M |\mathcal{N}_1(T,x^j_T) - g(x^j_T)|^2, \label{eq.loss2} \\
  \mbox{Loss}^3_{t_n} = \frac{1}{N}\frac{1}{M} \sum_{j=1}^M |\nabla\mathcal{N}_1(T,x^j_T) - \nabla g(x^j_T)|^2. \label{eq.loss3}
\end{gather}
It should be noted that in the temporal difference learning, which do not require the simulation of the entire trajectory, the loss function is computed at each time step. Hence, the terminal state $X_T$ remains unknown when $t<T$. To address this challenge, the terminal state $X_T$ from the previous iteration is utilized in the calculation of \eqref{eq.loss2} and \eqref{eq.loss3}, while for the initial iteration, the terminal state is randomly generated.

Lastly, as discussed in Subsection \ref{sec.neural}, equation \eqref{eq.N2loss} is incorporated as the fourth component of the loss function. 
\begin{equation} \label{eq.loss4}
  \begin{aligned}
    & \mbox{Loss}^4_{t_n} = \\
    & \left|
      \frac{1}{M} \sum_{j=1}^M\left(\sum_{i=N_n+1}^{N_{n+1}} [\mathcal{N}_1(t_n, x^j_{t_n}+G(x^j_{t_n},z^j_i))-\mathcal{N}_1(t_n,x^j_{t_n})] - \Delta t \mathcal{N}_2(t_n, x^j_{t_n})\right)
      \right|.
  \end{aligned}
\end{equation}
It serves to constrain the output of the neural network $\mathcal{N}_2$, aiming to make it as close to $\int_{\R^d}(\mathcal{N}_1(t,x+G(x,z)) - \mathcal{N}_1(t,x))\nu(\dif z)$ as possible. Based on these, the loss function at time step $t_n$ can be defined in the following manner
\begin{equation} \label{eq.loss}
  \mbox{Loss}_{t_n} = \mbox{Loss}^1_{t_n} + \mbox{Loss}^2_{t_n} + \mbox{Loss}^3_{t_n} + \mbox{Loss}^4_{t_n}.
\end{equation}

\begin{algorithm}
  \renewcommand{\algorithmicensure}{\textbf{Input:}}
	\renewcommand{\algorithmicrequire}{\textbf{Output:}}
  \caption{Temporal difference learning on high-dimensional PIDEs.} \label{algo.all}
  \begin{algorithmic}[1]
    \ENSURE{PIDE \eqref{eq.PIDE}, neural network with 2 outputs $ (\mathcal{N}_1,\mathcal{N}_2) $ and parameter $\theta$ , snapshots $ 0=t_0<t_1<t_2<\cdots<t_N=T $, sample number $ M $, maximum iterations $ Iterations $, initial distribution $ X_0=\xi $, learning rate $\alpha$.}
    \REQUIRE{Approximate solution $ \mathcal{N}_1 $ of PIDE \eqref{eq.PIDE}.}
    \STATE Initialize neural network. $iteration=0$. Randomly generate $M$ terminal state $\{x^j_T\}_{j=1}^M$ for the initial iteration. Initialize the previous terminal state $\mathbb{X}=\{x^j_T\}_{j=1}^M$.
    \WHILE {$iteration\leq Iterations$}
    \STATE Sample $M$ Brownian motion $\{\{\Delta W_n^j\}_{n=1}^N\}_{j=1}^M$ for each time interval $[t_n, t_{n+1}]$. Sample some jumps $\{z_i^j\}_{j=1}^M$ for each trajectory $j$ from L\'{e}vy motion (See section \ref{sec.numerical} for details).
    \FOR{time step $n=0\to N-1$}
    \STATE Update state $(t_n,x^j_n)$ to $(t_{n+1},x^j_{n+1})$ by \eqref{eq.dXt}, $j=1,2,\cdots,M$.
    \STATE Compute $\mbox{Loss}^1_{t_n}$ by \eqref{eq.loss1} and $\mbox{Loss}^4_{t_n}$ by \eqref{eq.loss4}.
    \IF{$n=N-1$}
    \STATE Update the previous terminal state $\mathbb{X}=\{x^j_T\}_{j=1}^M$.
    \ENDIF
    \STATE Compute $\mbox{Loss}^2_{t_n}$ by \eqref{eq.loss2} and $\mbox{Loss}^3_{t_n}$ by \eqref{eq.loss3} using the current $\mathbb{X}$.
    \STATE Compute $\mbox{Loss}_{t_n}$ by \eqref{eq.loss}. Optimize the parameter $\theta$ of the neural network.
    \ENDFOR
    \STATE $iteration=iteration+1$
    \ENDWHILE
  \end{algorithmic}
\end{algorithm}

\begin{figure}[htbp]
  \centering
  \includegraphics[width=37em]{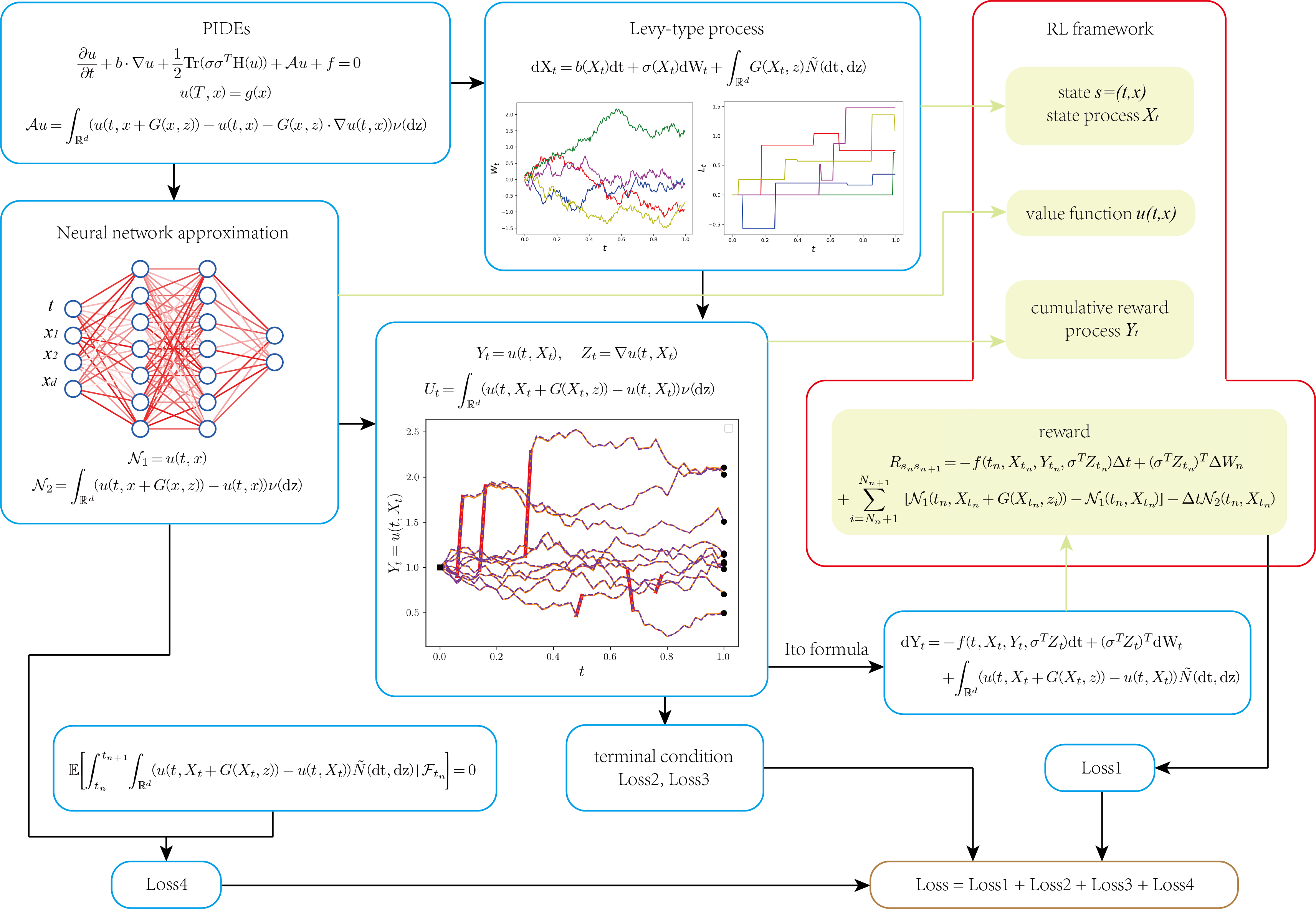}
  \caption{The diagram of solving high-dimensional PIDEs with jumps by temporal difference method. The PIDE under consideration is associated with a set of L\'{e}vy-type forward-backward stochastic processes. These processes can be effectively characterized and addressed within the framework of reinforcement learning, as depicted in the highlighted red box in the figure. By employing neural networks, the calculation of the loss function becomes feasible through the utilization of temporal difference methods. Consequently, the model parameters can be updated accordingly to facilitate the resolution of the PIDE.} \label{fig.diagram}
\end{figure}

The pseudo-code for the entire method is presented in Algorithm \ref{algo.all}. The overall process of this approach is illustrated in Figure \ref{fig.diagram}. In summary, the solution to the PIDE \eqref{eq.PIDE} is sought by introducing a set of L\'{e}vy-type forward-backward stochastic processes $(X_t, Y_t, Z_t, U_t)$. The simulation of $Y_t$ involves unknown functions $u(t,x)$ and $\int_{\R^d}(u(t,x+G(x,z)) - u(t,x))\nu(\dif z)$, which are represented separately using a residual network with two outputs. A reinforcement learning framework is then established, where $X_t$ corresponds to the state process, the solution $u(t,x)$ of \eqref{eq.PIDE} corresponds to the value function, $Y_t$ corresponds to the cumulative reward process, and the reward is obtained from equation \eqref{eq.dYt}. By employing temporal difference learning, termination conditions, and leveraging the properties of martingales, the final loss function is constructed. Finally, the parameters of the neural network are updated using optimization methods such as Adam. Through these steps, the solution to equation \eqref{eq.PIDE} is successfully obtained.

\begin{remark}
  From an alternative perspective, equation \eqref{eq.TD} employs a one-step temporal difference learning, where after the state transitions from $s_n$ to $s_{n+1}$, the sum of the reward $R_{s_ns_{n+1}}$ and $u(s_n)$ is used as an approximation for $u(s_{n+1})$. However, it is also possible to employ a two-step temporal difference learning
  \begin{equation}
    u(s_{n}) \leftarrow u(s_{n}) + \alpha(R_{s_ns_{n+1}} +R_{s_{n+1}s_{n+2}} + u(s_{n}) - u(s_{n+2})),
  \end{equation}
  whereby the state progresses from $s_n$ to $s_{n+2}$, and the sum of $R_{s_ns_{n+1}}+R_{s_{n+1}s_{n+2}}$ and $u(s_{n})$ is utilized to approximate $u(s_{n+2})$. Consequently, the TD error and loss1 would undergo corresponding adjustments. If, during a single iteration where the samples run from $t_0$ to $t_N$, the one-step temporal difference learning optimizes the parameters $2m$ times, then the two-step temporal difference learning optimizes the parameters $m$ times since it updates them every two steps. Similarly, this approach can be generalized to a $k$-step temporal difference learning.
\end{remark}


\section{Numerical results} \label{sec.numerical}

\begin{figure}[htbp]
  \centering
  \includegraphics[width=37em]{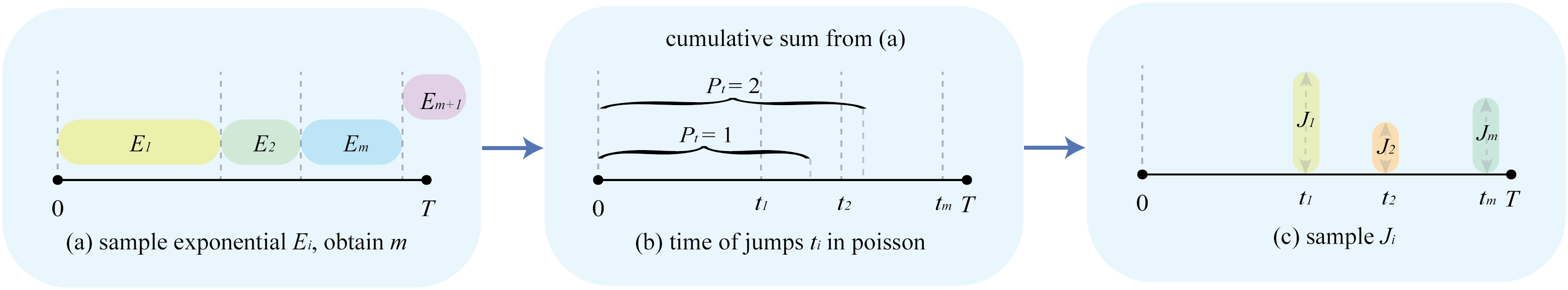}
  \caption{Simulation of jumps. To simulate the jumps on each trajectory, the following steps are taken: (a) A sequence of exponential distributions with parameter $\lambda$ is generated. (b) The cumulative sum of this sequence yields the arrival times of the Poisson process on the respective trajectory. (c) At each arrival time, a single sample is drawn from the distribution $\phi(z)$ to determine the value of the jump at that particular moment.} \label{fig.jumps}
\end{figure}

To facilitate the generation of the jump sizes $z_i$ in equations \eqref{eq.dXt} and \eqref{eq.dYt}, we consider the L\'{e}vy measure $\nu(\dif z)$ in the form of $\nu(\dif z)=\lambda\phi(z)\dif z$. Here $\phi(z)$ denotes the density function of a $d$-dimensional random variable and $\lambda\in\R$. By constructing a compound Poisson process $L_t$ using
\begin{equation}
  L_t = \sum_{i=1}^{P_t} J_i,
\end{equation}
where $\{J_i\}_{i=1}^\infty$ denotes a sequence of independently and identically distributed (i.i.d.) random variables with density function $\phi(z)$, and $P_t$ denotes a Poisson process with intensity $\lambda$ that is independent of $\{J_i\}_{i=1}^\infty$, the resulting jump measure $\nu(\dif z)$ associated with $L_t$ precisely corresponds to $\lambda\phi(z)\dif z$. 
Accordingly, the jumps can be simulated on each trajectory based on the procedure illustrated in Figure \ref{fig.jumps}. Initially, a sequence of exponential distributions $\{E_i\}_{i=1}^{m+1}$ with parameter $\lambda$ is generated. Subsequently, the cumulative sum of these random variables yields the arrival times $t_1,t_2,\cdots,t_m$ of the Poisson process $P_t$. Finally, the density function $\phi(z)$ is sampled $m$ times to simulate the jumps at each time point. This approach enables the determination of the jump occurrences times and their corresponding magnitudes.

In the numerical experiments, the neural network architecture is composed of a $d+1$ dimensional input, a linear layer, 5 residual blocks, a linear layer, and a 2-dimensional output, as illustrated in Figure \ref{fig.net}. Each residual block includes two linear layers mapping from 25-dimensional space to 25-dimensional space and a residual connection. In high-dimensional experiments, the width of the network is adjusted proportionally to the dimension. The \emph{tanh} activation function is employed, and the Adam optimizer is chosen with an initial learning rate of $5\times10^{-5}$. The learning rate is divided by 5 every 5000 parameter updates during the training process. All experiments adopt the same random seed of 2023. Considering the characteristics of backward stochastic differential equations, the $L_1$ relative error of $Y_0$ is employed as the evaluation of the model. The experiments are conducted on a Tesla V100 GPU with 32GB memory and 40GB RAM, utilizing 10 cores.

\subsection{One dimensional pure jump process} \label{sec.1d}

One-dimensional L\'{e}vy processes are widely encountered in various fields such as insurance, hedging, and pricing\cite{delong2013backward}. To better illustrate the abilities of our method in computing the PIDEs, this subsection specifically focuses on the one-dimensional pure jump process to avoid potential challenges arising from higher dimensions and Brownian motion. We aim to solve the following PIDE
\begin{equation}
    \left\{
    \begin{aligned}
        & \frac{\partial u}{\partial t}(t,x) + \int_{\R} \left(u(t,xe^z) - u(t,x) - x(e^z-1)\frac{\partial u}{\partial x}(t,x)\right)\nu(\dif z) = 0, \\ 
        & u(T, x) = x,
    \end{aligned}
    \right.
\end{equation}
where $ \nu(\dif z)=\lambda\phi(z)\dif z $ , $ \phi(z)=\frac{1}{\sqrt{2\pi}\sigma}e^{-\frac{1}{2}\left(\frac{z-\mu}{\sigma}\right)^2} $ and $x\in\R$ . The equation possesses an exact solution of $ u(t,x)=x $ , and the corresponding L\'{e}vy-type processes are
\begin{gather}
    \dif X_t = \int_{\R} X_t(e^z-1)\widetilde{N}(\dif t,\dif z), \\ 
    \dif Y_t = \int_{\R} \left(u(t,X_te^z) - u(t,X_t)\right)\widetilde{N}(\dif t,\dif z).
\end{gather}

The initial value of the process is set to $X_0=1$, and the time interval $[0,1]$ is divided into $N=50$ equal partitions. The intensity of the Poisson process is $\lambda=0.3$, while the jumps following a normal distribution with $\mu=0.4, \sigma=0.25$. A total of $M=1000$ trajectories are employed. The training consists of 400 iterations, with 50 parameter updates in each iteration, resulting in a cumulative total of 20,000 parameter updates.

\begin{figure}[h]
    \centering
    \subfigure[]{\includegraphics[height=13.5em]{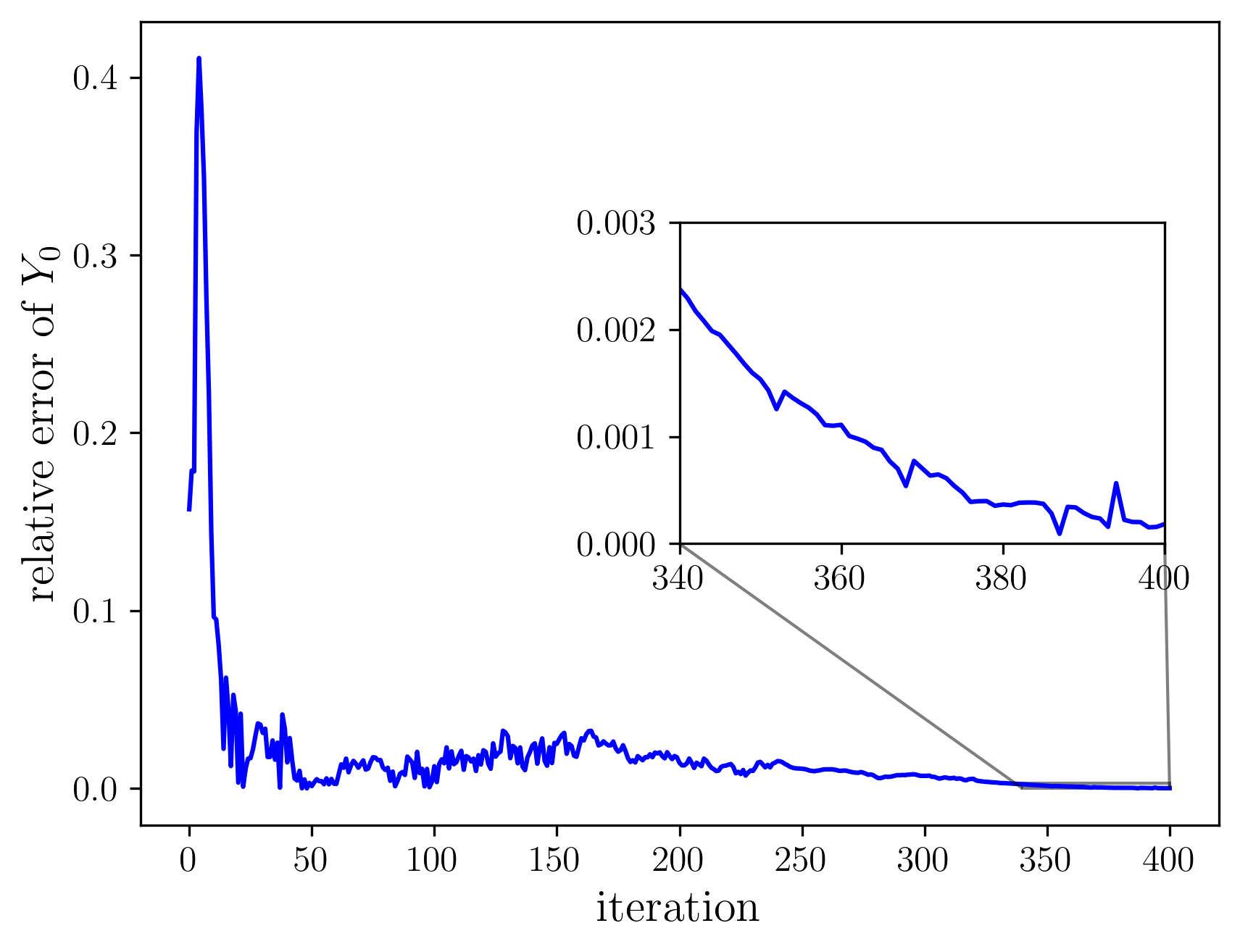}}
    \subfigure[]{\includegraphics[height=13.5em]{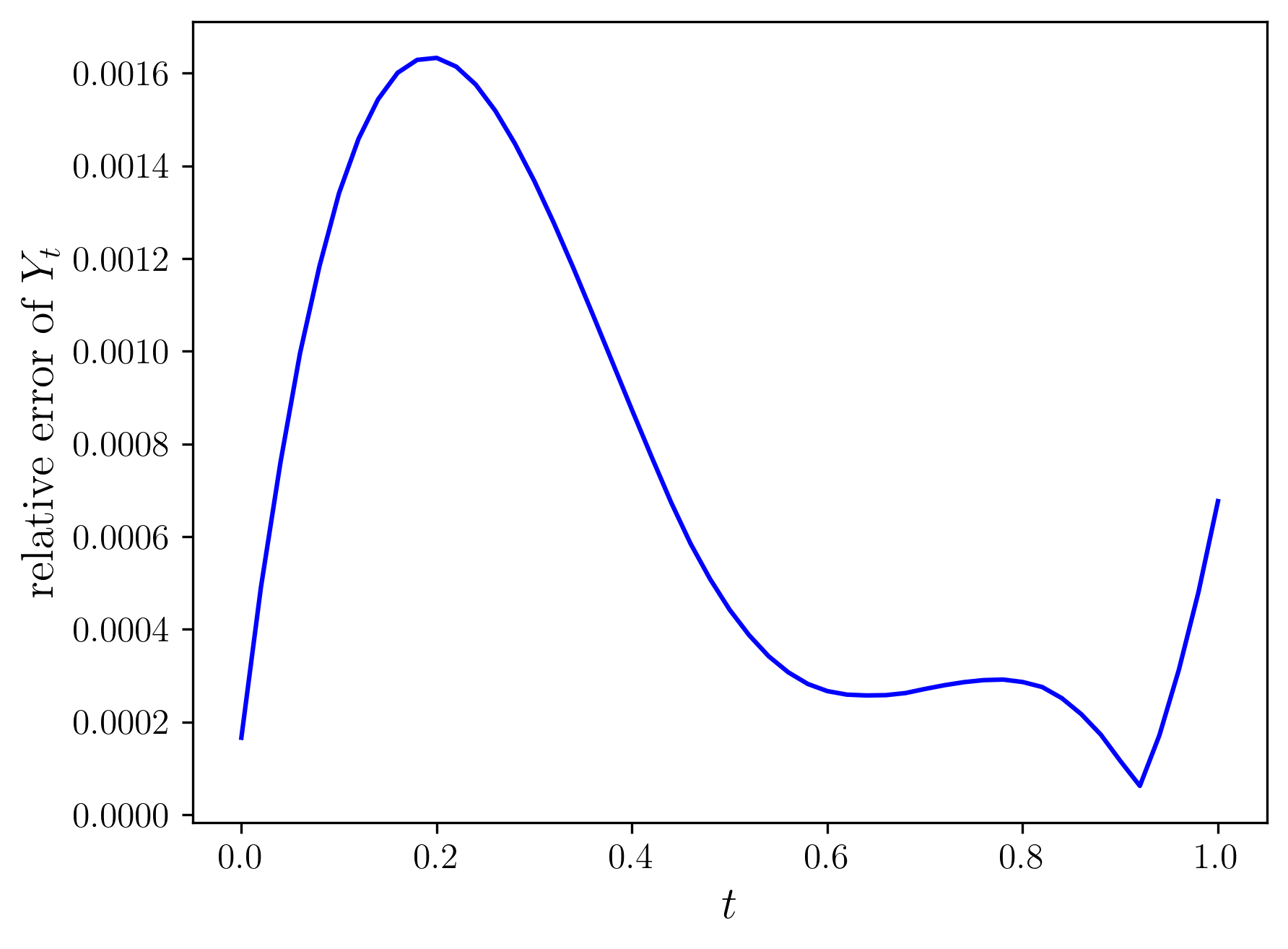}}
    \caption{Relative error of one-dimensional pure jump problem. (a) The evolution of the relative error of $Y_0$ with respect to the number of iterations. The exact value is $Y_0=1$, and the relative error converges to 0.02\%. (b) The relative error of the neural network's approximation for $Y_t$ at different time $t$.} \label{fig.1D_error}
\end{figure}

\begin{figure}[h]
    \centering
    \includegraphics[width=22em]{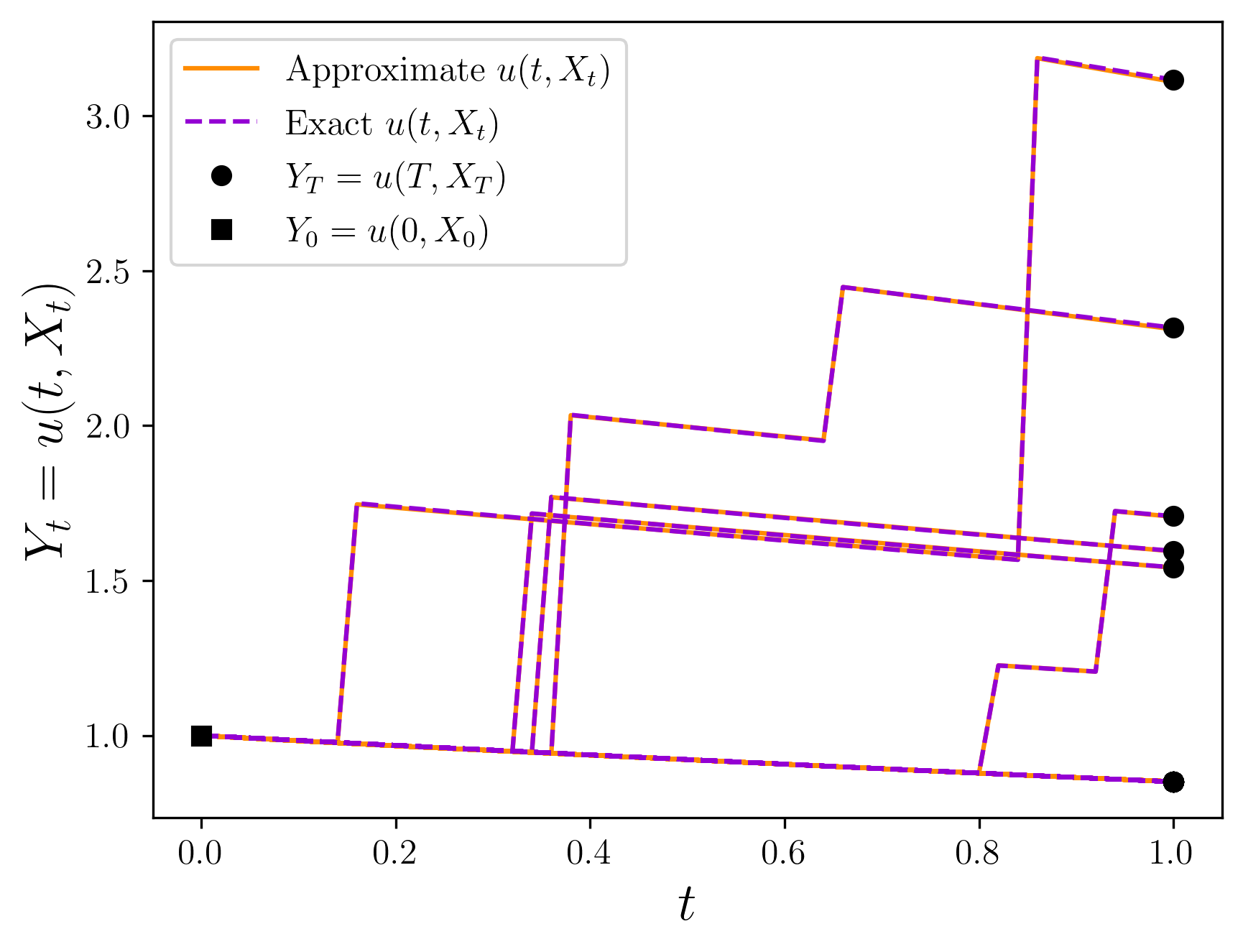}
    \caption{The visualization of the trajectories for the one-dimensional pure jump problem. 30 trajectories are displayed, and 5 trajectories exhibit jumps, while the remaining 25 trajectories coincide entirely due to the absence of Brownian motion.} \label{fig.1D_traj}
\end{figure}

Figure \ref{fig.1D_error}(a) illustrates the evolution of the relative error of $Y_0$ throughout the entire training process, while Table \ref{tab.1D} presents the specific values of the relative error of $Y_0$ during the initial and final 5000 steps of training. Figure \ref{fig.1D_error}(b) depicts the relative error of the neural network's approximation for $Y_t$ at each time $t$. Our method achieves a relative error of 0.02\% for $Y_0$. It is noteworthy that the error decreases to approximately 1\% after only 1000 parameter updates, and with sufficient training, the error magnitude reaches the order of $10^{-4}$. To visually assess the results, Figure \ref{fig.1D_traj} illustrates 30 trajectories, among which 5 exhibit jumps. The results illustrate the rapid convergence and high precision achieved by our method in one-dimensional pure jump problems. In summary, our method demonstrates remarkable performance.

\begin{table}[htbp]
    \centering
    \caption{Relative error of $Y_0$ during the initial and final 5000 updates of training in the one-dimensional pure jump problem. The exact value is $Y_0=1$.} \label{tab.1D}
    \begin{tabular}{ccc|ccc}
        \toprule
        Iteration & $Y_0$ & Relative error & Iteration & $Y_0$ & Relative error \\ 
        \midrule
        0 & 1.1570 & 15.70\% & 15000 & 0.9929 & 0.71\% \\
        500 & 1.0965 & 9.65\% & 15500 & 0.9943 & 0.57\% \\ 
        1000 & 1.0033 & 0.33\% & 16000 & 0.9946 & 0.54\% \\ 
        1500 & 1.0313 & 3.13\% & 16500 & 0.9968 & 0.32\% \\ 
        2000 & 0.9853 & 1.47\% & 17000 & 0.9976 & 0.24\% \\ 
        2500 & 1.0013 & 0.13\% & 17500 & 0.9985 & 0.15\% \\ 
        3000 & 0.9974 & 0.26\% & 18000 & 0.9989 & 0.11\% \\ 
        3500 & 0.9862 & 1.38\% & 18500 & 0.9993 & 0.07\% \\ 
        4000 & 0.9896 & 1.04\% & 19000 & 0.9996 & 0.04\% \\ 
        4500 & 0.9837 & 1.63\% & 19500 & 0.9997 & 0.03\% \\ 
        5000 & 0.9875 & 1.25\% & 20000 & 0.9998 & 0.02\% \\
        \bottomrule
    \end{tabular}
\end{table}

\subsection{Robustness}

The proposed method involves various parameters, which may not always align with the ideal conditions encountered in practical scenarios. The effectiveness of our method has been demonstrated in one-dimensional pure jump problem. In this subsection, the performance of the method under different parameter settings is evaluated to assess its robustness, including the number of trajectories, the number of time intervals, jump intensity, jump form, TD step size. The experiment is conducted using the following PIDE:
\begin{equation} \label{eq.robust}
    \left\{
    \begin{aligned}
        \frac{\partial u}{\partial t}(t,x) + \epsilon & x\frac{\partial u}{\partial x}(t,x) + \frac{1}{2}\theta^2\frac{\partial^{2} u}{\partial x^{2}} (t,x) \\
        & + \int_{\R} \left(u(t,xe^z) - u(t,x) - x(e^z-1)\frac{\partial u}{\partial x}(t,x)\right)\nu(\dif z) = \epsilon x, \\ 
        u(T, x) = x,
    \end{aligned}
    \right.
\end{equation}
where $\nu(\dif z)=\lambda\phi(z)\dif z$, and $x\in\R$ . The equation possesses an exact solution of $ u(t,x)=x $, and the corresponding L\'{e}vy-type processes are
\begin{gather}
    \dif X_t = \epsilon X_t \dif t + \theta\dif W_t + \int_{\R} X_t(e^z-1)\widetilde{N}(\dif t,\dif z), \\ 
    \dif Y_t = \epsilon X_t \dif t + \theta Z_t\dif W_t + \int_{\R} \left(u(t,X_te^z) - u(t,X_t)\right)\widetilde{N}(\dif t,\dif z),
\end{gather}
and the initial value is set to $X_0=1$.

\textbf{(a) Trajectories and intervals.} 
In the computation of the loss function, Monte Carlo methods are employed to estimate the expectations. Therefore, the number of trajectories $M$ used can potentially influence the method's performance. The number of time intervals $N$ is associated with the discretization of the stochastic process, as well as the number of optimization each iteration. Hence, it is an important factor in assessing the method's robustness. Here we assume the jump form to follow a normal distribution $\phi(z)=\frac{1}{\sqrt{2\pi}\sigma}e^{-\frac{1}{2}\left(\frac{z-\mu}{\sigma}\right)^2}$ with $\mu=0.4, \sigma=0.25$ and the intensity of the Poisson process is set to $\lambda=0.3$, maintaining consistency with subsection \ref{sec.1d}. Additionally, we set $\epsilon=0.25, \theta=0$ in equation \eqref{eq.robust} and the termination time $T=1$.

Since the number of time intervals $N$ affects the number of optimization performed in each iteration, we report the results for three different numbers of iterations: 250, 500, and 750. For example, in the case of 500 iterations and $N=40$ time intervals, this implies a total of $500*40=20,000$ parameter updates.
The results are presented in Table \ref{tab.robust_traj}. The findings reveal that increasing the number of trajectories $M$ can effectively reduce the error at convergence, with a more pronounced improvement observed when the number of trajectories increases from $M=125$ to $M=250$. On the other hand, altering the number of time intervals between $N=20, 40, 80$ does not exert a significant influence on the error at convergence. For trajectory numbers of $M=250$ or higher, the error at convergence almost consistently remains within the order of $10^{-3}$. Therefore, our method does not require a high number of trajectories and is not sensitive to the number of time intervals.

\begin{table}[htbp]
    \centering
    \caption{Relative errors of $Y_0$ for different numbers of trajectories $M$ and time intervals $N$ in the robustness checks. Parameters are updated $N$ times in each iteration.} \label{tab.robust_traj}
    {\fontsize{6}{8}\selectfont
    \begin{tabular}{c|ccc|ccc|ccc}
        \toprule
        iterations & \multicolumn{3}{c}{250} & \multicolumn{3}{|c|}{500} & \multicolumn{3}{c}{750} \\ 
        \midrule
        intervals $N$ & 20 & 40 & 80 & 20 & 40 & 80 & 20 & 40 & 80 \\  
        \midrule
        $M=1000$ & 7.687\% & 1.671\% & 1.723\% & 0.303\% & 0.660\% & 0.545\% & 0.276\% & 0.159\% & 0.377\% \\
        $M=500$ & 2.988\% & 0.412\% & 2.611\% & 0.245\% & 1.385\% & 0.665\% & 0.219\% & 0.402\% & 0.639\% \\ 
        $M=250$ & 0.972\% & 1.671\% & 2.847\% & 0.264\% & 0.598\% & 0.502\% & 2.122\% & 0.139\% & 0.364\% \\ 
        $M=125$ & 4.100\% & 5.687\% & 4.303\% & 3.272\% & 0.662\% & 1.968\% & 1.700\% & 0.489\% & 1.691\% \\ 
        \bottomrule
    \end{tabular}
    }
\end{table}


\textbf{(b) L\'{e}vy measure.} 
The form of the L\'{e}vy measure plays a decisive role in non-local terms. Here, we aim to evaluate the performance of our method under different Poisson intensities $\lambda$ and jump form $\phi(z)$. A total of $M=250$ trajectories are employed, and the time interval $[0,1]$ is evenly divided into $N=50$ intervals. The parameters in Equation \eqref{eq.robust} are set as $\epsilon=0, \theta=0.4$. The training process comprises 400 iterations, resulting in a total of 20,000 parameter updates. For the conducted tests, the Poisson intensities are systematically varied from 0.3 to 1.8. Additionally, the following four distinct jump distributions $\phi(z)$ are utilized: 1) normal distribution with $\mu=0.4, \sigma=0.25$, 2) uniform distribution with $\delta=0.4$, 3) exponential distribution with $\lambda_0=3$, and 4) Bernoulli distribution with $a_1=-0.2, a_2=0.4, p=0.7$:
\begin{equation} \label{eq.jump_form}
    \begin{gathered}
        \phi(z)=\frac{1}{\sqrt{2\pi}\sigma}e^{-\frac{1}{2}\left(\frac{z-\mu}{\sigma}\right)^2}, \quad
        \phi(z)=\begin{cases} \frac{1}{2\delta}, & -\delta\leq z\leq\delta \\ 0, & \mbox{else}\end{cases}, \\ 
        \phi(z)=\begin{cases} \lambda_0 e^{-\lambda_0z}, & z\geq0 \\ 0, & z<0\end{cases}, \quad 
        \phi(z)=\begin{cases} p & z=a_1 \\ 1-p, & z=a_2\end{cases}.
    \end{gathered}
\end{equation}

\begin{figure}[htbp]
    \centering
    \subfigure[Normal distribution]{\includegraphics[width=18em]{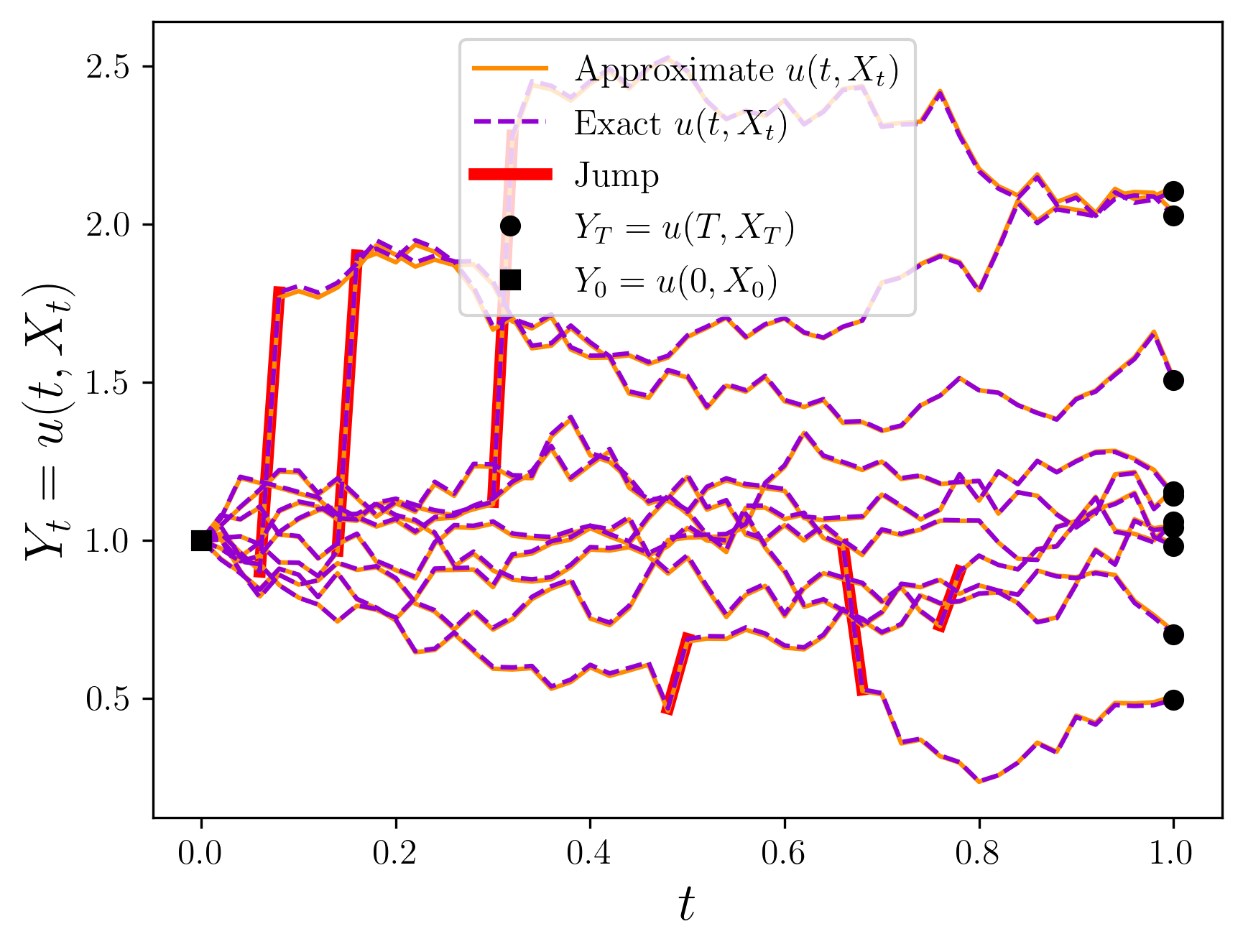}}
    \subfigure[Uniform distribution]{\includegraphics[width=18em]{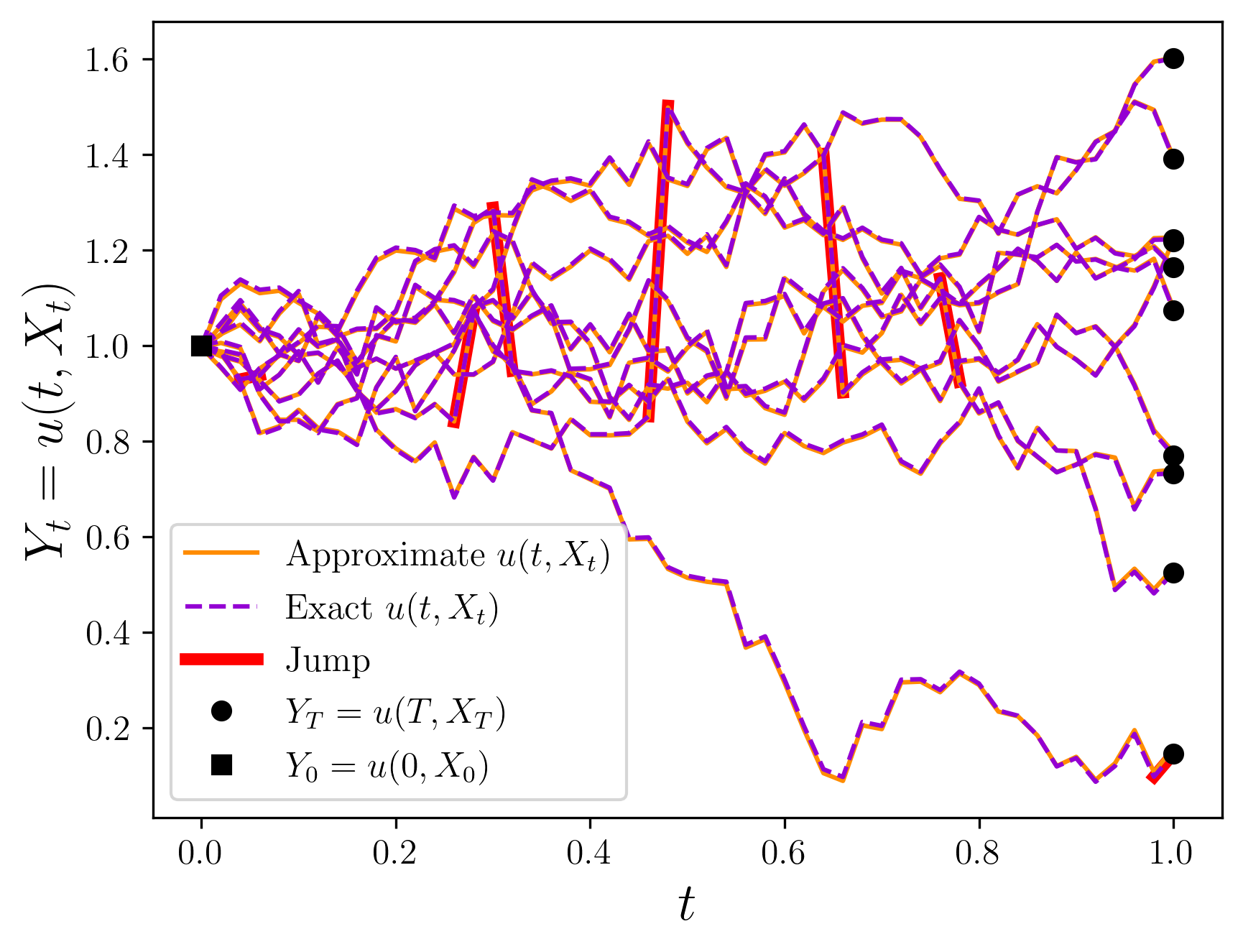}}

    \subfigure[Exponential distribution]{\includegraphics[width=18em]{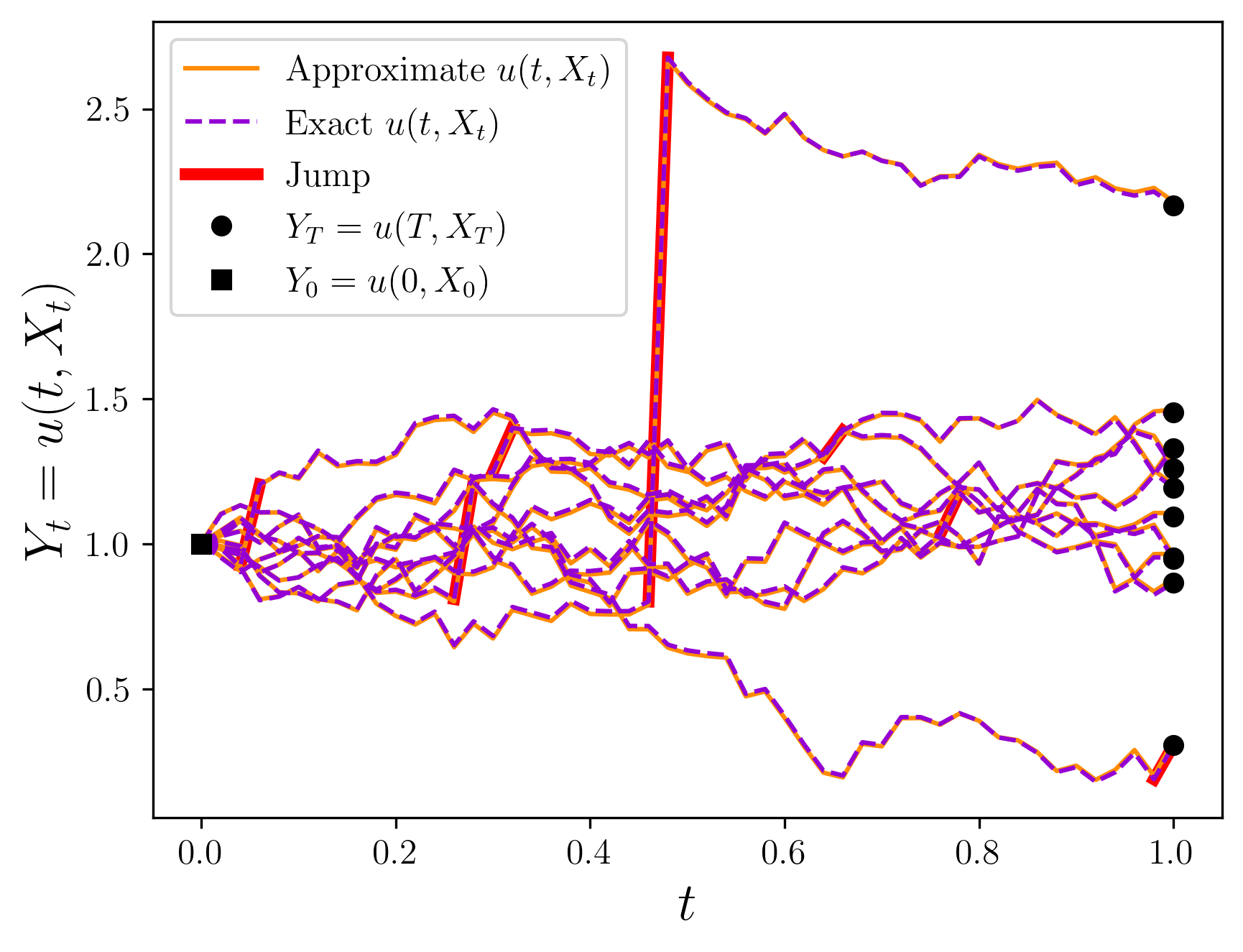}}
    \subfigure[Bernoulli distribution]{\includegraphics[width=18em]{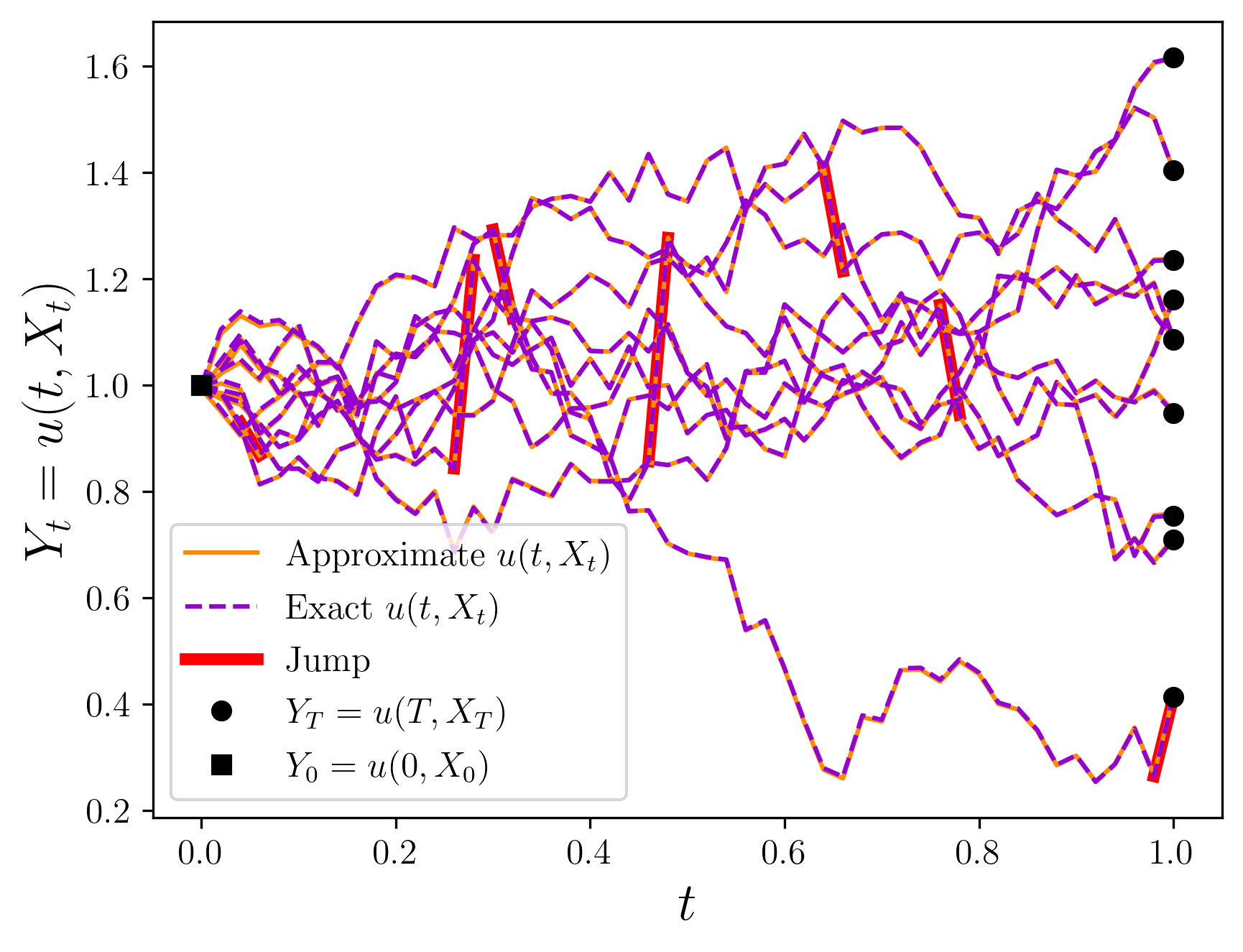}}
    \caption{The visualization of the trajectories under different jump forms $\phi(z)$ in the robustness checks. To distinguish the jumps with smaller amplitudes from the Brownian motion, the occurrence of jumps is marked with bold red lines.}
    \label{fig.robust_jumps}
\end{figure}

Table \ref{tab.robust_jumps} illustrates the relative errors of $Y_0$ under various Poisson intensities $\lambda$ and jump forms $\phi(z)$. With the increase of Poisson intensity from $\lambda=0.3$ to $\lambda=1.8$, the occurrence frequency of jumps within the $[0, T]$ time interval gradually increases. The results reveal a growing trend in the error for the normal distribution, while the errors for the remaining three distributions show no significant changes. Across almost all scenarios, the relative errors of $Y_0$ remain within 2\%. Figure \ref{fig.robust_jumps} visually showcases 10 trajectories for each of the four jump forms when the Poisson intensity is $\lambda=0.3$. Notably, the occurrence of jumps is indicated by bold red lines. The figure illustrates the excellent performance of our method under the four distinct jump forms.

\begin{table}[htbp]
    \centering
    \caption{Relative errors of $Y_0$ under different Poisson intensities $\lambda$ and jump forms $\phi(z)$ in the robustness checks. The specific forms of the four distributions are given in Equation \eqref{eq.jump_form}. The experiments were all conducted with $M=250$ trajectories, $N=50$ time intervals, and a total of 20,000 parameter updates.} \label{tab.robust_jumps}
    \begin{tabular}{ccccc}
        \toprule
        $Y_0$ relative error & Normal & Exponential & Uniform & Bernoulli \\
        \midrule
        $\lambda=0.3$ & 0.051\% & 0.475\% & 0.376\% & 0.752\% \\
        $\lambda=0.6$ & 0.117\% & 0.393\% & 2.450\% & 2.469\% \\ 
        $\lambda=0.9$ & 0.340\% & 0.910\% & 1.355\% & 1.997\% \\ 
        $\lambda=1.2$ & 1.111\% & 0.567\% & 1.048\% & 1.202\% \\ 
        $\lambda=1.5$ & 1.960\% & 1.591\% & 1.202\% & 1.227\% \\ 
        $\lambda=1.8$ & 1.173\% & 0.461\% & 0.167\% & 1.257\% \\ 
        \bottomrule
    \end{tabular}
\end{table}

\textbf{(c) $n$-step Temporal difference.} 
According to the loss function, the previous experiments can be regarded as one-step TD methods. However, an intuitive approach is to investigate the influence of selecting an appropriate TD step $n$ on the method's performance. In this regard, we set $\epsilon=0.25, \theta=0.4$ in Equation \eqref{eq.robust} for training purposes. A total of $M=250$ trajectories were employed, and the Poisson process intensity was set to $\lambda=0.3$, with the jump form being $\phi(z)=\frac{1}{\sqrt{2\pi}\sigma}e^{-\frac{1}{2}\left(\frac{z-\mu}{\sigma}\right)^2}$ and $\mu=0.4, \sigma=0.25$. For the sake of computational convenience, it was desirable to have the number of intervals $N$ precisely divisible by the TD step $n$. Consequently, the time $[0,1]$ was divided into $N=60$ intervals, and we examined the TD step $n=1, 2, 3, 4, 5, 6$.

Due to the reduced number of parameter updates each iteration with larger step sizes $n$, and considering both the number of parameter updates and convergence, we performed $M=500$ iterations for $n=1,2$, $M=1000$ iterations for $n=3,4$, and $M=1500$ iterations for $n=5,6$. The results of the relative error of $Y_0$ for different step sizes $n$ are presented in Table \ref{tab.robust_TD}. It can be observed that, in this case, the one-step TD method demonstrates the most favorable computational performance. The error does not exhibit a significant trend as the step size increases, and it remains almost consistently below 1\%.

\begin{table}[htbp]
    \centering
    \caption{Relative errors of $Y_0$ for different TD step $n$ in the robustness checks. We performed $M=500$ iterations for $n=1,2$, $M=1000$ iterations for $n=3,4$, and $M=1500$ iterations for $n=5,6$.} \label{tab.robust_TD}
    \begin{tabular}{ccccccc}
        \toprule
        TD step $n$ & 1 & 2 & 3 & 4 & 5 & 6 \\ 
        \midrule
        $Y_0$ & 0.9998 & 0.9926 & 1.0114 & 1.0088 & 1.0078 & 0.9965 \\
        relative error of $Y_0$ & 0.022\% & 0.745\% & 1.137\% & 0.880\% & 0.780\% & 0.354\% \\
        \bottomrule
    \end{tabular}
\end{table}

\subsection{High dimensional problems}
When faced with higher-dimensional problems, many methods encounter a bottleneck due to the exponential increase in computational cost. This phenomenon, known as the curse of dimensionality, poses a significant challenge. In this subsection, our focus shifts towards investigating the performance in high-dimensional scenarios. The experiments are conducted using the following PIDE
\begin{equation} \label{eq.highdim}
  \left\{
  \begin{aligned}
      & \frac{\partial u}{\partial t}(t,x) + \frac{\epsilon}{2}x\cdot\nabla u(t,x) + \frac{1}{2}\mbox{Tr}(\theta^2\mbox{H}(u)) \\
      & + \int_{\R^d} \left(u(t,x+z) - u(t,x) - z\cdot\nabla u(t,x)\right)\nu(\dif z)
      = \lambda(\mu^2+\sigma_0^2) + \theta^2 + \frac{\epsilon}{d}\|x\|^2, \\ 
      & u(T, x) = \frac{1}{d}\|x\|^2, \\ 
  \end{aligned}
  \right.
\end{equation}
with $\epsilon=0, \theta=0.3, x\in\R^d$ and Poisson intensity $\lambda=0.3$. To highlight the effectiveness of our method in handling high-dimensional problems, we chose the jump form to be independent multi-dimensional normal distributions, with a mean of $\mu*(1,1,\cdots,1)$ and a covariance matrix $\sigma^2 I_d$. Here $\mu=0.1, \sigma_0^2=0.0001$, and the intention is to maintain uniformity in jump sizes to the greatest extent possible. The equation possesses an exact solution of $u(t, x) = \frac{1}{d}\|x\|^2$. The training process utilizes $M=500$ trajectories, and the time $[0, 1]$ is divided into $N=50$ intervals. Furthermore, the number of neurons in the linear layers of each block is adjusted to $d+10$, instead of 25, to accommodate the varying input dimensions $d$.

\begin{figure}[htbp]
  \centering
  \subfigure[]{\includegraphics[width=18.2em, height=13.6em]{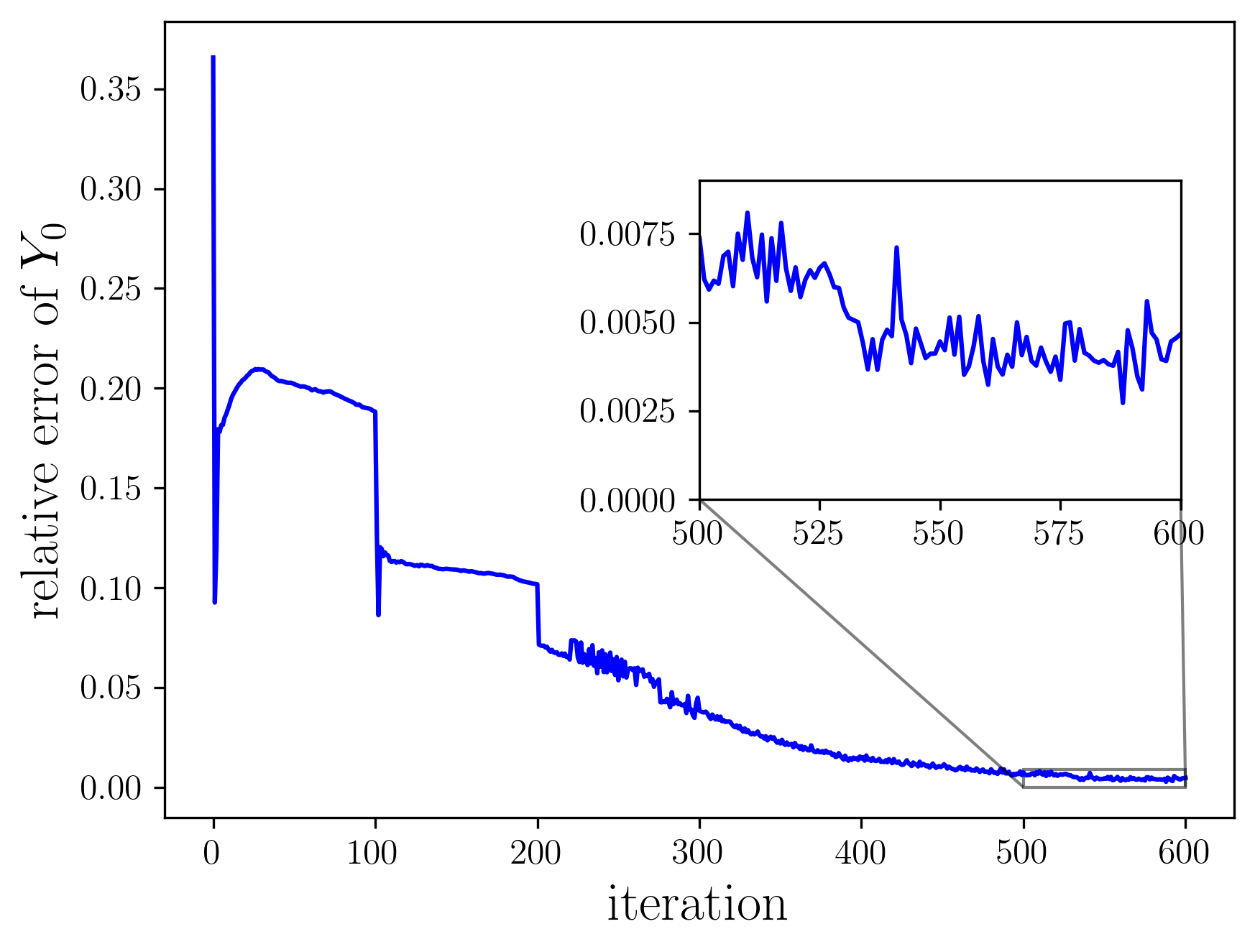}}
  \subfigure[]{\includegraphics[width=17.8em, height=13.8em]{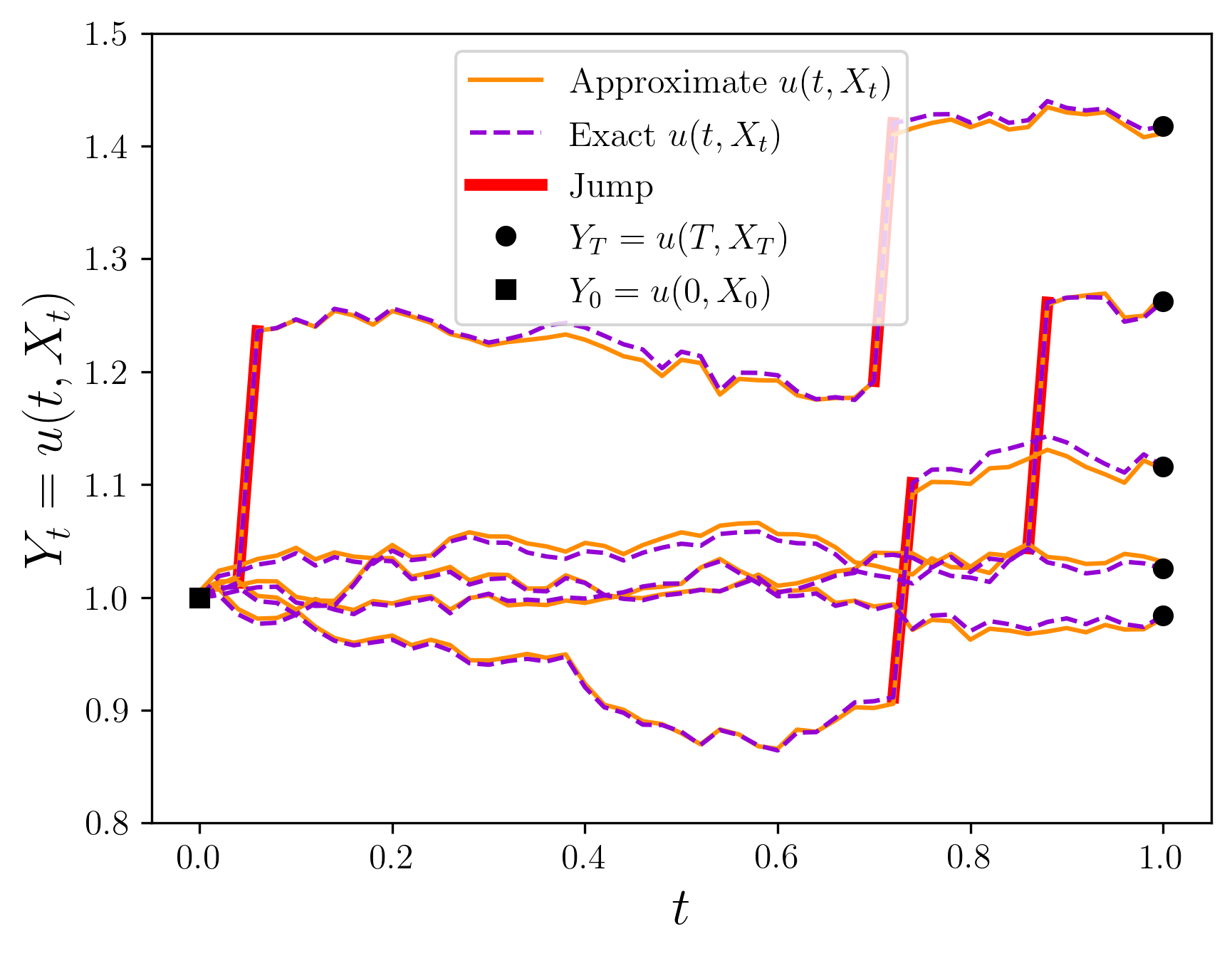}}
  \caption{Result of the 100-dimensional problem. The exact value of $Y_0$ is 1, and the relative error of the neural network approximation is 0.548\%. (a) The evolution of relative error of $Y_0$ with the number of iterations. (b) Visualization of 5 trajectories under the neural network's approximation and the exact solution.} \label{fig.100d}
\end{figure}

In the 100-dimensional experiment, the relative error of $Y_0$ is found to be 0.548\% with a computational time of approximately 10 minutes. Figure \ref{fig.100d}(a) depicts the evolution of the relative error of $Y_0$ throughout the
entire training process in the 100-dimensional experiment. Since the model did not converge after 400 iterations, i.e., 20,000 parameter updates, the training was extended to 600 iterations, i.e., 30,000 parameter updates. Furthermore, Figure \ref{fig.100d}(b) provides a visual representation of 5 trajectories in the 100-dimensional problem, highlighting the locations of jumps through bold red lines. It is evident that our proposed method achieves a high-precision approximation of the solution in the 100-dimensional problem, while maintaining a computationally feasible runtime.

Furthermore, the performance of our method was evaluated in various dimensions $d$. The experimental setup remained consistent with the 100-dimensional experiment. Table \ref{tab.highdim} presents the relative error of $Y_0$ and the program's execution time across different dimensions $d$. It is noteworthy that regardless of low or high dimensions, the relative error of $Y_0$ generally remains within the range of $10^{-3}$. All experiments were conducted using $M=500$ trajectories, and it is evident that our method exhibits a minimal sensitivity to the number of trajectories in high-dimensional problems. As the dimensionality increases, the program's execution time exhibits a roughly linear growth, ranging from 5 to 10 minutes, which is entirely acceptable. In summary, our method demonstrates remarkable performance across different dimensions and shows no significant curse of dimensionality.

\begin{table}[htbp]
  \centering
  \caption{The relative error of $Y_0$ and runtime for different dimensions $d$. The experiments were conducted on a Tesla V100 GPU with 32GB memory
  and 40GB RAM.} \label{tab.highdim}
  \begin{tabular}{cccccc}
      \toprule
      Dimension & 2 & 4 & 6 & 8 & 10 \\
      $Y_0$ relative error & 0.954\% & 0.251\% & 0.025\% & 0.671\% & 1.895\% \\
      Time (s) & 332 & 336 & 338 & 356 & 364 \\
      \midrule
      Dimension & 20 & 30 & 40 & 50 & 60 \\
      $Y_0$ relative error & 0.702\% & 1.221\% & 0.956\% & 0.219\% & 0.944\% \\
      Times (s) & 368 & 396 & 401 & 423 & 506 \\
      \midrule
      Dimension & 70 & 80 & 90 & 100 & - \\
      $Y_0$ relative error & 0.044\% & 0.277\% & 0.460\% & 0.548\% & - \\
      Time (s) & 529 & 565 & 584 & 638 & - \\
      \bottomrule
  \end{tabular}
\end{table}

\begin{figure}[htbp]
  \centering
  \subfigure[$d=25$]{\includegraphics[width=18em]{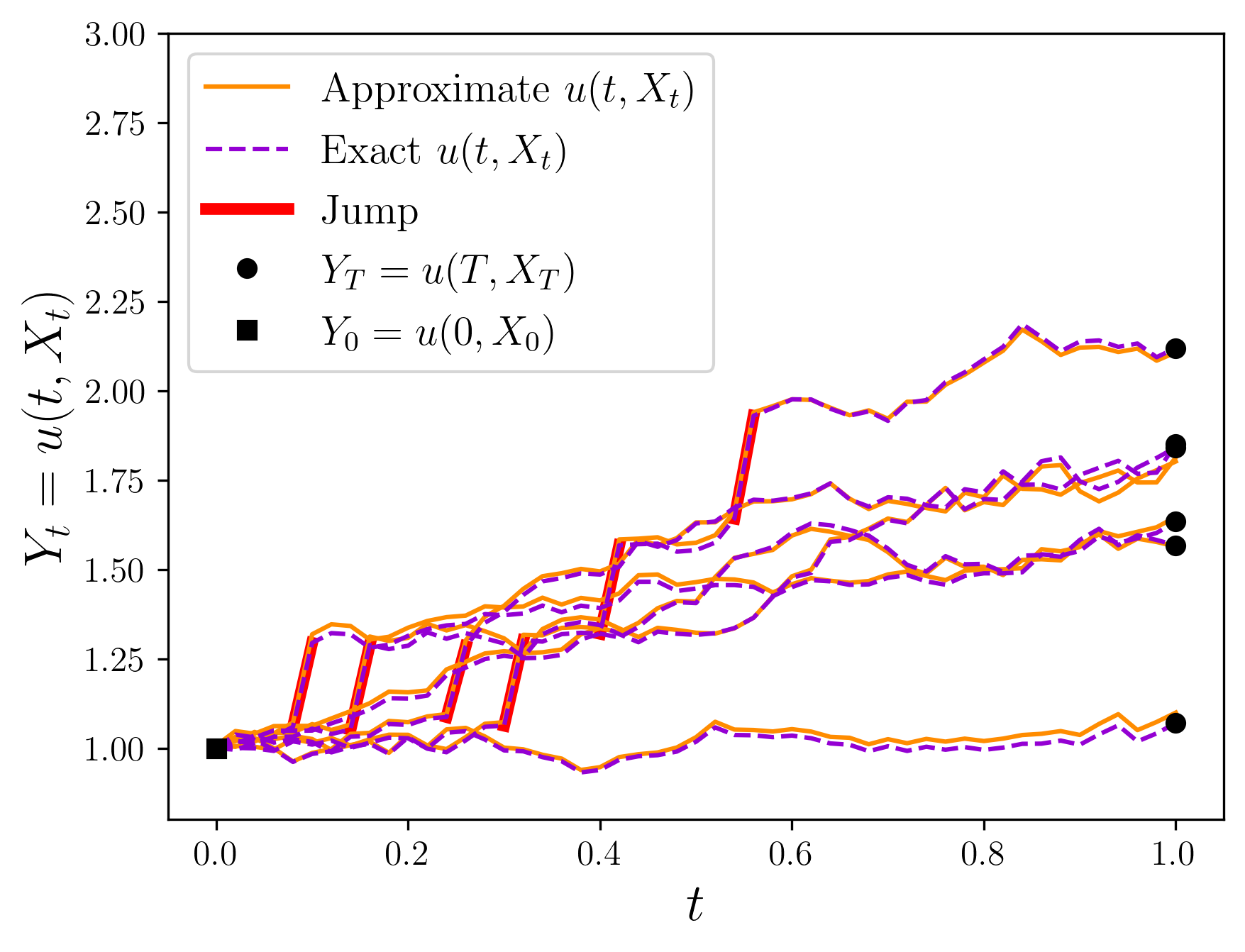}}
  \subfigure[$d=50$]{\includegraphics[width=18em]{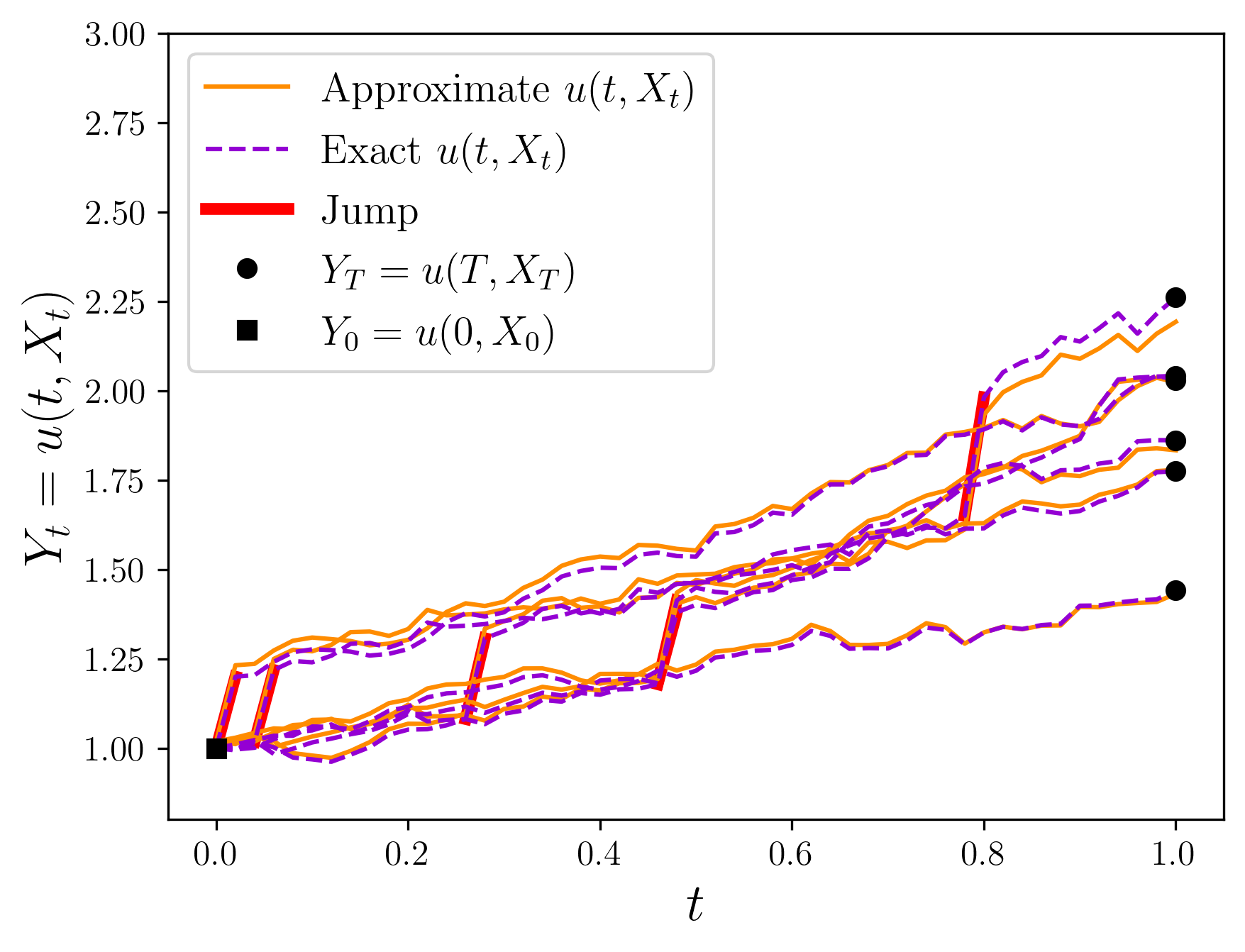}}

  \subfigure[$d=75$]{\includegraphics[width=18em]{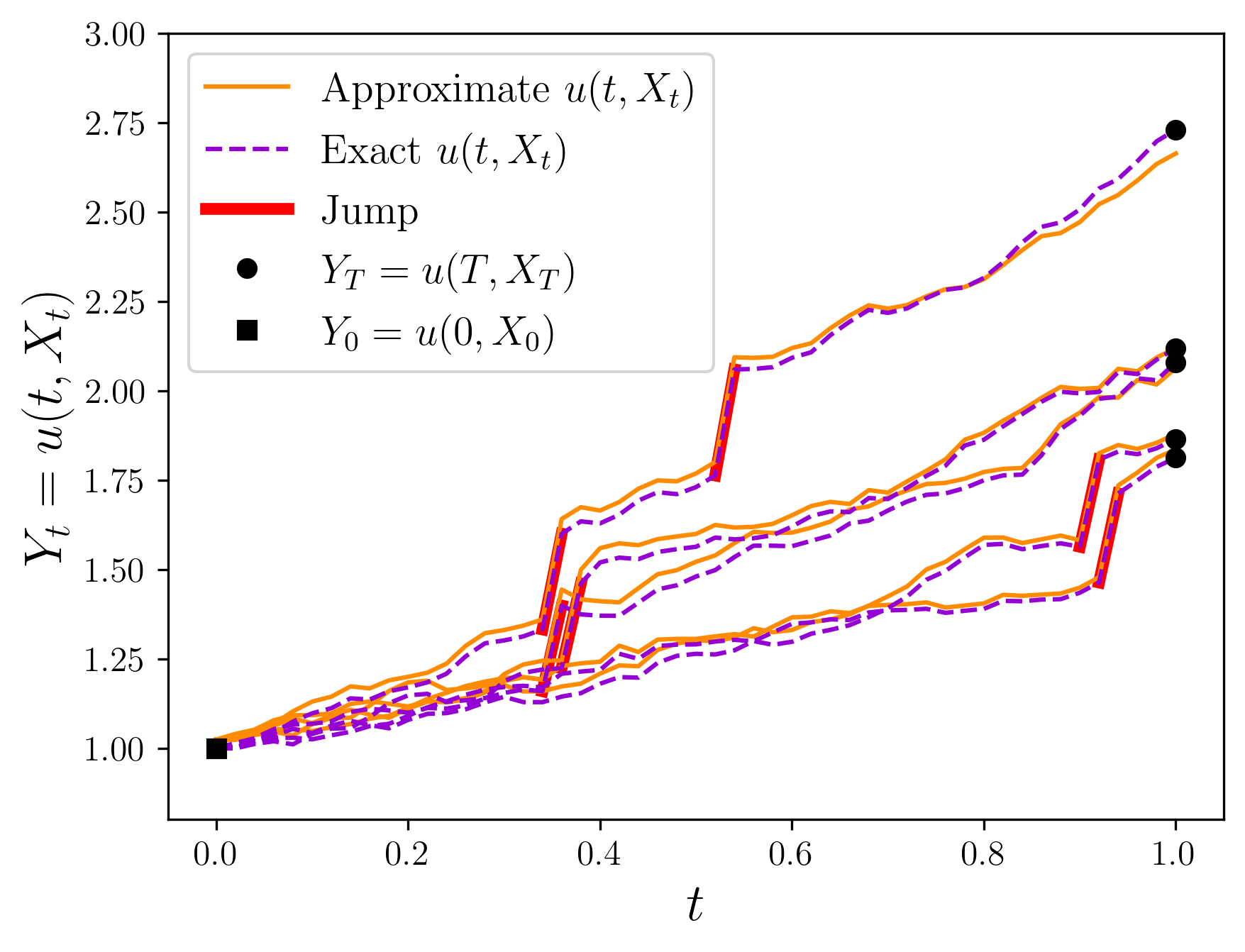}}
  \subfigure[$d=100$]{\includegraphics[width=18em]{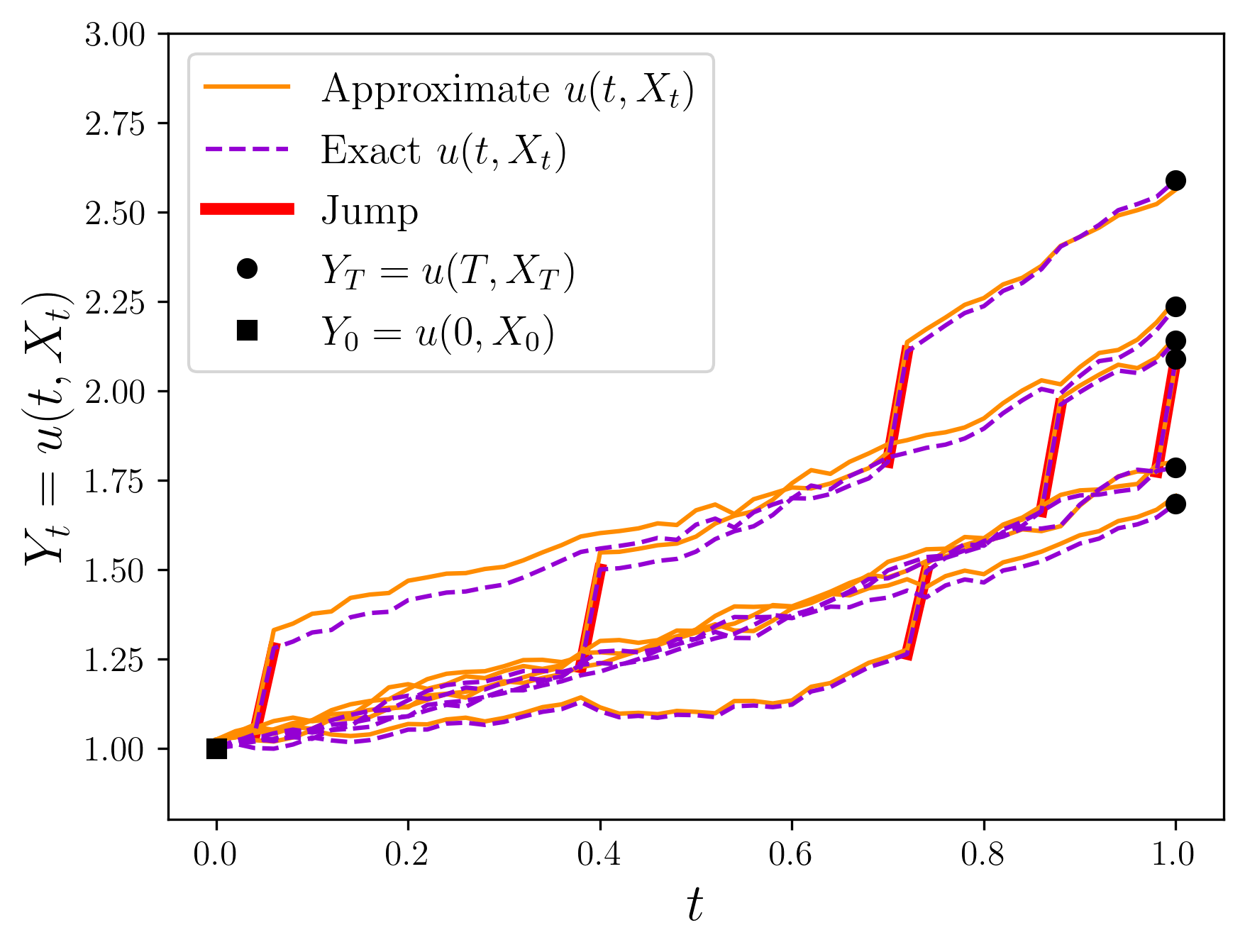}}
  \caption{The visualization of the trajectories of high dimensional problems with coupled drift term and couple diffusion term.}
  \label{fig.highdim_couple}
\end{figure}

Finally, we slightly modified equation \eqref{eq.highdim} to consider the following high dimensional PIDE with coupled drift term and coupled diffusion term
\begin{small}
\begin{equation}
  \left\{
  \begin{aligned}
      & \frac{\partial u}{\partial t}(t,x) + \frac{\epsilon}{2}\left \| x \right \| x\cdot\nabla u(t,x) + \frac{1}{2}\mbox{Tr}(\sigma\sigma^T\mbox{H}(u)) \\
      & + \int_{\R^d} \left(u(t,x+z) - u(t,x) - z\cdot\nabla u(t,x)\right)\nu(\dif z)
      = \lambda(\mu^2+\sigma_0^2) + \frac{2d-2}{d}\theta^2 + \frac{\epsilon}{d}\|x\|^3,  \\ 
      & u(T, x) = \frac{1}{d}\|x\|^2. \\ 
  \end{aligned}
  \right.
\end{equation}
\end{small}
Here, $\epsilon=0.05$, $\theta=0.2$, Poisson intensity $\lambda=0.3$, and the diffusion term is given by
\[
  \sigma = \theta
  \begin{pmatrix}
    1 & 0 & 0 & 0 & \cdots & 0 \\
    1 & 1 & 0 & 0 & \cdots & 0 \\
    0 & 1 & 1 & 0 & \cdots & 0 \\
    0 & 0 & 1 & 1 & \cdots & 0 \\
    \vdots & \vdots & \vdots & \vdots & \ddots & 0 \\ 
    0 & 0 & 0 & \cdots & 1 & 1 
  \end{pmatrix}.
\]
The exact solution is $u(t,x)=\frac{1}{d}\|x\|^2$. In this numerical test, we use the same setups as in the previous example, including the selection of the jump form, the discretization of the time interval, the number of trajectories, and the number of iterations. We test the cases for dimensions $d=25, \;50,\; 75$, and $100$, and visualize their trajectories in Figure \ref{fig.highdim_couple}. The relative errors of $Y_0$ are 0.743\% for $d=25$, 1.910\% for $d=50$, 2.412\% for $d=75$ and 2.387\% for $d=100$. It is seen that our method still performs well with the existence of the coupled drift term and coupled diffusion term.

\section{Conclusion} \label{sec.conclusion}

In this work, a framework based on reinforcement learning is proposed to solve high-dimensional partial integro-diﬀerential equations. A set of L\'{e}vy-type stochastic processes is introduced to characterize the solutions of the PIDE, and a reinforcement learning model is built upon these processes. Deep neural networks are employed to represent the solutions of the original PIDE and non-local terms, enabling a step-by-step simulation of the system. The temporal difference learning is utilized, and its error is employed as the loss function, allowing parameter updates without the need to wait for completing the entire trajectory. Additionally, termination conditions and the properties of non-local terms are incorporated into the loss function. The entire training process can be completed with a relatively low computational cost.

The numerical experiments demonstrate the efficacy of the proposed method in solving 100-dimensional problems with an relative error on the order of $10^{-3}$ and one-dimensional pure jump problems with an relative error on the order of $10^{-4}$. The computational cost of the method exhibits approximately linear growth with increasing dimensionality, while maintaining a low requirement for the number of trajectories in high-dimensional cases. Robustness tests reveal favorable outcomes for different forms and intensities of jumps, without the need for dense time interval divisions or large temporal difference step sizes. These findings indicate the broader applicability of the proposed method.
The error analysis is left for the future work. Furthermore, the extension of this method to more general problems, such as inverse problems involving L\'{e}vy processes with available data, warrants further investigation in future work.

\section*{Acknowledgements}
L.L. and Y.Z. were partially supported by the National Key R\&D Program of China (Grant No. 2021YFA0719200). H.G. was partially supported by the Andrew Sisson Fund, Dyason Fellowship, and the Faculty Science Researcher Development Grant of the University of Melbourne. X.Y. was partially supported by the NSF grant DMS-2109116. L.L. would like to thank Sheng Wang for helpful discussions.



\bibliographystyle{siamplain}
\bibliography{References}

\begin{thebibliography}{10}

\bibitem{abergel2010nonlinear}
{\sc F.~Abergel and R.~Tachet}, {\em A nonlinear partial integro-differential equation from mathematical finance}, Discrete and Continuous Dynamical Systems-Series A, 27 (2010), pp.~907--917.

\bibitem{applebaum2009levy}
{\sc D.~Applebaum}, {\em L{\'e}vy processes and stochastic calculus}, Cambridge university press, 2009.

\bibitem{cai2021physics}
{\sc S.~Cai, Z.~Mao, Z.~Wang, M.~Yin, and G.~E. Karniadakis}, {\em Physics-informed neural networks (pinns) for fluid mechanics: A review}, Acta Mechanica Sinica, 37 (2021), pp.~1727--1738.

\bibitem{castro2022deep}
{\sc J.~Castro}, {\em Deep learning schemes for parabolic nonlocal integro-differential equations}, Partial Differential Equations and Applications, 3 (2022), p.~77.

\bibitem{chiu2022can}
{\sc P.-H. Chiu, J.~C. Wong, C.~Ooi, M.~H. Dao, and Y.-S. Ong}, {\em {CAN-PINN}: A fast physics-informed neural network based on coupled-automatic--numerical differentiation method}, Computer Methods in Applied Mechanics and Engineering, 395 (2022), p.~114909.

\bibitem{cont2005finite}
{\sc R.~Cont and E.~Voltchkova}, {\em A finite difference scheme for option pricing in jump diffusion and exponential {L}{\'e}vy models}, SIAM Journal on Numerical Analysis, 43 (2005), pp.~1596--1626.

\bibitem{cuomo2022scientific}
{\sc S.~Cuomo, V.~S. Di~Cola, F.~Giampaolo, G.~Rozza, M.~Raissi, and F.~Piccialli}, {\em Scientific machine learning through physics--informed neural networks: Where we are and what’s next}, Journal of Scientific Computing, 92 (2022), p.~88.

\bibitem{delong2013backward}
{\sc {\L}.~Delong}, {\em Backward stochastic differential equations with jumps and their actuarial and financial applications}, Springer, 2013.

\bibitem{duan2015introduction}
{\sc J.~Duan}, {\em An introduction to stochastic dynamics}, vol.~51, Cambridge University Press, 2015.

\bibitem{frey2022deep}
{\sc R.~Frey and V.~K{\"o}ck}, {\em Deep neural network algorithms for parabolic {PIDEs} and applications in insurance mathematics}, in Mathematical and Statistical Methods for Actuarial Sciences and Finance: MAF 2022, Springer, 2022, pp.~272--277.

\bibitem{gnoatto2022deep}
{\sc A.~Gnoatto, M.~Patacca, and A.~Picarelli}, {\em A deep solver for {BSDEs} with jumps}, arXiv:2211.04349,  (2022).

\bibitem{goswami2016system}
{\sc A.~Goswami, J.~Patel, and P.~Shevgaonkar}, {\em A system of non-local parabolic {PDE} and application to option pricing}, Stochastic Analysis and Applications, 34 (2016), pp.~893--905.

\bibitem{goswami2013optimal}
{\sc D.~Goswami, A.~K. Pani, and S.~Yadav}, {\em Optimal error estimates of two mixed finite element methods for parabolic integro-differential equations with nonsmooth initial data}, Journal of Scientific Computing, 56 (2013), pp.~131--164.

\bibitem{han2018solving}
{\sc J.~Han, A.~Jentzen, and W.~E}, {\em Solving high-dimensional partial differential equations using deep learning}, Proceedings of the National Academy of Sciences, 115 (2018), pp.~8505--8510.

\bibitem{han2017deep}
{\sc J.~Han, A.~Jentzen, et~al.}, {\em Deep learning-based numerical methods for high-dimensional parabolic partial differential equations and backward stochastic differential equations}, Communications in mathematics and statistics, 5 (2017), pp.~349--380.

\bibitem{karniadakis2021physics}
{\sc G.~E. Karniadakis, I.~G. Kevrekidis, L.~Lu, P.~Perdikaris, S.~Wang, and L.~Yang}, {\em Physics-informed machine learning}, Nature Reviews Physics, 3 (2021), pp.~422--440.

\bibitem{kavallaris2018non}
{\sc N.~I. Kavallaris and T.~Suzuki}, {\em Non-local partial differential equations for engineering and biology}, Mathematical Modeling and Analysis, 31 (2018).

\bibitem{kharazmi2021hp}
{\sc E.~Kharazmi, Z.~Zhang, and G.~E. Karniadakis}, {\em hp-{VPINNs}: Variational physics-informed neural networks with domain decomposition}, Computer Methods in Applied Mechanics and Engineering, 374 (2021), p.~113547.

\bibitem{koppen2000curse}
{\sc M.~K{\"o}ppen}, {\em The curse of dimensionality}, in 5th online world conference on soft computing in industrial applications (WSC5), vol.~1, 2000, pp.~4--8.

\bibitem{kwon2011second}
{\sc Y.~Kwon and Y.~Lee}, {\em A second-order finite difference method for option pricing under jump-diffusion models}, SIAM journal on numerical analysis, 49 (2011), pp.~2598--2617.

\bibitem{li2022revisiting}
{\sc W.~Li, C.~Zhang, C.~Wang, H.~Guan, and D.~Tao}, {\em Revisiting {PINNs}: Generative adversarial physics-informed neural networks and point-weighting method}, arXiv:2205.08754,  (2022).

\bibitem{liang2022finite}
{\sc S.~Liang and H.~Yang}, {\em {Finite Expression Method for Solving High-Dimensional Partial Differential Equations}}, arXiv:2206.10121,  (2022).

\bibitem{liao2019deep}
{\sc Y.~Liao and P.~Ming}, {\em Deep nitsche method: Deep ritz method with essential boundary conditions}, Communications in Computational Physics, 29 (2021), pp.~1365--1384.

\bibitem{lu2021deepxde}
{\sc L.~Lu, X.~Meng, Z.~Mao, and G.~E. Karniadakis}, {\em {DeepXDE}: A deep learning library for solving differential equations}, SIAM review, 63 (2021), pp.~208--228.

\bibitem{lyu2022mim}
{\sc L.~Lyu, Z.~Zhang, M.~Chen, and J.~Chen}, {\em {MIM}: A deep mixed residual method for solving high-order partial differential equations}, Journal of Computational Physics, 452 (2022), p.~110930.

\bibitem{mcclenny2023self}
{\sc L.~D. McClenny and U.~M. Braga-Neto}, {\em Self-adaptive physics-informed neural networks}, Journal of Computational Physics, 474 (2023), p.~111722.

\bibitem{moseley2021finite}
{\sc B.~Moseley, A.~Markham, and T.~Nissen-Meyer}, {\em {Finite Basis Physics-Informed Neural Networks (FBPINNs): a scalable domain decomposition approach for solving differential equations}}, arXiv:2107.07871,  (2021).

\bibitem{pablo2023finite}
{\sc O.~Pablo and D.~Ciro}, {\em A finite elements approach for spread contract valuation via associated two-dimensional {PIDE}}, Computational and Applied Mathematics, 42 (2023), p.~15.

\bibitem{raissi2018forward}
{\sc M.~Raissi}, {\em Forward-backward stochastic neural networks: Deep learning of high-dimensional partial differential equations}, arXiv:1804.07010,  (2018).

\bibitem{raissi2019physics}
{\sc M.~Raissi, P.~Perdikaris, and G.~E. Karniadakis}, {\em Physics-informed neural networks: A deep learning framework for solving forward and inverse problems involving nonlinear partial differential equations}, Journal of Computational physics, 378 (2019), pp.~686--707.

\bibitem{sirignano2018dgm}
{\sc J.~Sirignano and K.~Spiliopoulos}, {\em {DGM}: A deep learning algorithm for solving partial differential equations}, Journal of computational physics, 375 (2018), pp.~1339--1364.

\bibitem{wang2023deep}
{\sc W.~Wang, J.~Wang, J.~Li, F.~Gao, and Y.~Fu}, {\em Deep learning numerical methods for high-dimensional fully nonlinear {PIDEs} and coupled {FBSDEs} with jumps}, arXiv:2301.12895,  (2023).

\bibitem{wight2020solving}
{\sc C.~L. Wight and J.~Zhao}, {\em Solving allen-cahn and cahn-hilliard equations using the adaptive physics informed neural networks}, Communications in Computational Physics, 29 (2021), pp.~930--954.

\bibitem{yu2018deep}
{\sc B.~Yu et~al.}, {\em The deep {R}itz method: a deep learning-based numerical algorithm for solving variational problems}, Communications in Mathematics and Statistics, 6 (2018), pp.~1--12.

\bibitem{yu2022gradient}
{\sc J.~Yu, L.~Lu, X.~Meng, and G.~E. Karniadakis}, {\em Gradient-enhanced physics-informed neural networks for forward and inverse {PDE} problems}, Computer Methods in Applied Mechanics and Engineering, 393 (2022), p.~114823.

\bibitem{zang2020weak}
{\sc Y.~Zang, G.~Bao, X.~Ye, and H.~Zhou}, {\em Weak adversarial networks for high-dimensional partial differential equations}, Journal of Computational Physics, 411 (2020), p.~109409.

\bibitem{zeng2022deep}
{\sc S.~Zeng, Y.~Cai, and Q.~Zou}, {\em Deep neural networks based temporal-difference methods for high-dimensional parabolic partial differential equations}, Journal of Computational Physics, 468 (2022), p.~111503.

\end{thebibliography}
\end{document}